\documentclass[journal]{IEEEtran}
\usepackage{amssymb}
\usepackage{amsmath}
\usepackage{amsfonts}
\usepackage{graphicx}
\usepackage{graphics}
\usepackage{subfigure}
\usepackage{epsfig}
\usepackage{stmaryrd}
\usepackage{mathrsfs}
\usepackage{latexsym}
\usepackage{float}
\usepackage{color}
\usepackage{cases}
\usepackage{cite}
\usepackage{enumerate}
\usepackage{times}

\newtheorem{theorem}{Theorem}
\newtheorem{corollary}{Corollary}
\newtheorem{definition}{Definition}
\newtheorem{lemma}{Lemma}

\newtheorem{assumption}{Assumption}
\newtheorem{remark}{Remark}
\newtheorem{example}{Example}

\ifCLASSINFOpdf
\else
\fi
\long\def\comment#1{}

\newcommand{\diag}{\mathrm{diag}}
\newcommand{\rank}{\mathrm{rank}}
\newcommand{\sgn}{\mathrm{sgn}}

\graphicspath{{figures/}}

\begin{document}

\title{Generalized Dynamics in Social Networks With Antagonistic Interactions}

\author{Shidong~Zhai~~and~~Wei~Xing~Zheng

%
\thanks{S.~Zhai is with the Chongqing Key Laboratory of Complex Systems and Bionic Control, Chongqing University of Posts and Telecommunications, Chongqing, 400065, China, and also with the School of Computing, Engineering and Mathematics, Western Sydney University, Sydney, NSW 2751, Australia (e-mail: zhaisd@cqupt.edu.cn).}
\thanks{W.~X.~Zheng is with the School of Computing, Engineering and Mathematics, Western Sydney University, Sydney, NSW 2751, Australia (e-mail: w.zheng@westernsydney.edu.au).}
}
\date{}

\maketitle

\begin{abstract}
In this paper, we investigate a general nonlinear model of opinion dynamics in which both state-dependent susceptibility to persuasion and antagonistic interactions are considered. According to the existing literature and socio-psychological theories, we examine three specializations of state-dependent susceptibility, that is, stubborn positives scenario, stubborn neutrals scenario, and stubborn extremists scenario. Interactions among agents form a signed graph, in which positive and negative edges represent friendly and antagonistic interactions, respectively. Based on Perron-Frobenius property of eventually positive matrices and LaSalle invariance principle, we conduct a comprehensive theoretical analysis of the generalized nonlinear opinion dynamics. We obtain some sufficient conditions such that the states of all agents converge into the subspace spanned by the right positive eigenvector of an eventually positive matrix. When there exists at least one entry of the right positive eigenvector which is not equal to one, the derived results can be used to describe different levels of an opinion. Finally, we present two examples to demonstrate the effectiveness of the theoretical findings.
\end{abstract}


\section{Introduction}

\IEEEPARstart{T}{he} study of opinion dynamics has attracted much attention in recent years. Based on some social psychology theories such as social comparison theory \cite{festinger1954theory}, cognitive dissonance theory \cite{festinger1962theory} and balance theory \cite{cartwright1956structural,heider1946attitudes}, various mathematical models have been presented to model the evolution of opinion forming \cite{degroot1974reaching,pirani2014spectral,friedkin1999social,friedkin2015problem,%
jia2015opinion,ghaderi2014opinion,groeber2014dissonance,xu2015modified}. Among these models, the basis is DeGroot model \cite{degroot1974reaching}, which is a discrete-time and weighted averaging model. The corresponding continuous-time version of DeGroot model was studied in \cite{Abelson64}. The system matrix of DeGroot model is the row-stochastic adjacency matrix of an interaction graph, and cannot change with the states. Moreover, DeGroot model cannot model social behavior of stubborn agents, which is often hard to be affected by other agents. In \cite{pirani2014spectral} DeGroot model, in which some agents' states are unchanged, was investigated.

To overcome the limitation of DeGroot model, Friedkin-Johnsen model \cite{friedkin1999social,friedkin2015problem} permits the agents to have different susceptibilities to persuasion. Although Friedkin-Johnsen model can model the behavior of stubborn agents, the behavior of agents cannot change with time and do not depend on current opinion or attitude. Various extensions of Friedkin-Johnsen model have been considered in \cite{jia2015opinion,ghaderi2014opinion,groeber2014dissonance,xu2015modified}.
A typical opinion model with state-dependent interactions is the bounded confidence model \cite{deffuant2000mixing,hegselmann2002opinion,lorenz2007continuous,castellano2009statistical,%
lorenz2005stabilization,blondel2009krause,mirtabatabaei2012opinion}.
This model assumes that two agents can interact with each other when they are close enough. Some special bounded confidence models have been studied in \cite{jalili2013effects,sobkowicz2015extremism,chen2016characteristics}.

Recently, Amelkin et al.~\cite{amelkin2017polar} proposed a continuous-time nonlinear model of polar opinion dynamics in which the agents' susceptibilities to persuasion depend on the agents' current opinion or attitude at hand, and interactions among agents are cooperative. According to different theories of social psychology, they investigated three specialized scenarios, that is, stubborn extremists, stubborn positives and stubborn neutrals. The scenario with stubborn extremists, based on social-psychological references \cite{miller1993attitude,eagly1993psychology,bassili1996meta}, assumes that the extreme opinions are resistant to change and neutral opinions will be attracted by extreme opinions, while the scenario with stubborn positives, based on \cite{friedkin2015problem}, is a special version of model with stubborn extremists. The scenario with stubborn neutrals, based on social comparison theory \cite{festinger1954theory} and social norms \cite{friedkin2015problem}, supposes that the extreme opinions are more likely to change than the neutral opinions. Using the non-smooth analysis tool, they made a thorough theoretical analysis of the novel opinion dynamics model in \cite{amelkin2017polar}.

In this paper, we will study a generalized nonlinear opinion dynamics which includes both state-dependent susceptibility to persuasion and antagonistic interactions. Cooperative and antagonistic interactions usually coexist in social networks \cite{wasserman1994social,easley2010networks,altafini2012dynamics}. The signed graphs are often used to describe the cooperative and antagonistic interactions in social networks, with positive and negative edges denoting respectively the cooperative and antagonistic interactions. In social networks, the signed graphs are generally structurally balanced \cite{cartwright1956structural} in the sense that all nodes of the graph can be divided into two subsets, where the edges in each subset are positive and the edges connecting the two subsets are negative. Based on the balance theory \cite{cartwright1956structural}, the bipartite consensus for structurally balanced signed graphs has been investigated in \cite{altafini2013consensus,HuZheng2013CDC,zhang2014bipartite,hu2014emergent,%
valcher2014consensus,zhai2016pinning,qin2017TCYB}.
Generally speaking, when the signed graph is structurally balanced, the bipartite consensus problem can be transformed to a classical consensus or synchronization problem. However, when the signed graph is not structurally balanced, one should not expect the network to achieve bipartite consensus. Instead, the unanimity of opinion over a signed graph was studied in \cite{altafini2015predictable} when the adjacency matrix satisfies the eventually positive property \cite{noutsos2006perron,olesky2009m,elhashash2008general}.

Our introduced model can be seen as a generalization of the opinion forming models in \cite{amelkin2017polar,altafini2015predictable}. Comparing it with the model in \cite{amelkin2017polar}, our model can deal with the antagonistic interactions and accurately describe the degree of agent's opinion. For example, ten people voted for the Green Party in an election, but they did not support it at the same level. Hence, it is reasonable to use different positive or negative values to differentiate the opinion of these electors. Comparing it with the model in \cite{altafini2015predictable}, our model introduces state-dependent susceptibility to persuasion. This modification makes the dynamics of our model become more complex than the one in \cite{altafini2015predictable}. There also exist some unobservable behaviors in the linear model in \cite{altafini2015predictable} (see Example~\ref{exm2} in Section~\ref{section4} for details).

This paper considers three specializations of state-dependent susceptibility which are drawn from \cite{amelkin2017polar}. We assume that all eigenvalues of the coupled system matrix have non-positive real parts, and the dominant eigenvalue is real and the corresponding eigenvector is a positive vector. We consider the directed and undirected signed graphs, respectively. Two cases about the adjacency matrix of a signed graph are examined: 1) the adjacency matrix is eventually positive; 2) the adjacency matrix is not eventually positive, but there exists a matrix such that the adjacency matrix has Perron-Frobenius property by adding it. When the dominant eigenvalue of the coupled system matrix is zero, utilizing Perron-Frobenius property of matrix with some negative entries and LaSalle invariance principle, we obtain sufficient conditions such that the states of the studied model will converge into the subspace spanned by the right positive eigenvector of the adjacency matrix. On the other hand, when the dominant eigenvalue of the coupled system matrix is non-positive, we derive sufficient conditions such that the general equilibrium points of the underlying model are asymptotically stable.

The remainder of the paper is organized as follows. Section~\ref{section2} presents notations and the problem statement. Section~\ref{section3} undertakes a comprehensive theoretical analysis of three specializations of state-dependent susceptibility, that is, stubborn positives scenario, stubborn neutrals scenario, and stubborn extremists scenario in generalized nonlinear opinion dynamics. Section~\ref{section4} gives two numerical examples to verify the obtained theoretical results. Finally, Section~\ref{section5} summarizes our conclusions and describes future work. Besides,  Appendix~\ref{appendix1} provides some definitions and facts about signed graphs and Appendix~\ref{appendix2} states some results of eventually positive matrices and Perron-Frobenius property, followed by some lemmas about positive semidefinite matrices given in Appendix~\ref{appendix3}.

\section{Preliminaries}\label{section2}

\subsection{Notations}

Let $\mathbf{1}\, (\mathbf{0})$ be a vector of appropriate dimension with all elements equal to $1 (0)$. Let $\mathbb{R}^{n}_{+}:=\{x\mid x\geq\mathbf{0}\}$, $\mathbb{R}^{n}_{-}:=\{x\mid x\leq\mathbf{0}\}$, and $\mathbb{R}^{n}_{-,+}:=\{x\mid x\leq\mathbf{0}\;\text{or}\;x\geq\mathbf{0}\}$, where $x=[x_1\;x_2\;\cdots\;x_n]^{\rm T}$. $\sgn$ is the signum function. For a vector $x\in\mathbb{R}^n$, let $\diag(x):=\diag\{x_1,x_2,\cdots,x_n\}$, $\diag(x)^2:=\diag\{x_1^2,x_2^2,\cdots,x_n^2\}$, and $||x||_{\infty}:=\max_{i}|x_i|$ which is the infinity norm of $x$. For vectors $x, y\in\mathbb{R}^n$, $x<y$ if $x_i<y_i, i=1,2,\cdots,n$. For a real square matrix $A=[a_{ij}]\in\mathbb{R}^{n\times n}$, $A\geq0$ if $a_{ij}\geq0$, and $A>0$ if $a_{ij}>0$, $i,j=1,2,\cdots,n$. If $A\geq0\, (>0)$, then $A$ is called a nonnegative (positive) matrix. For a real symmetric square matrix $B$, $B\succeq0\, (\preceq0)$ means that $B$ is positive semidefinite (negative semidefinite), and $B\succ0\, (\prec0)$ means that $B$ is positive definite (negative definite). $I$ represents an identity matrix of appropriate dimension. Let $\text{sp}(A)=\{\lambda_1(A),\cdots,\lambda_n(A)\}$ denote the spectrum of matrix $A$, and $\max\{\lambda_1(A),\cdots,\lambda_n(A)\}$ denote the dominant eigenvalue of matrix $A$, where $\lambda_i(A), i=1,\cdots,n$ are the eigenvalues of $A$. $\rho(A)$ stands for the spectral radius of matrix $A$, which is the smallest real positive number such that $\rho(A)\geq |\lambda_i(A)|, \forall i=1,\cdots,n$. Denote the transpose of matrix $A$ by $A^{\rm T}$, and let $H(A)=\frac{1}{2}(A+A^{\rm T})$. A matrix $A$ is irreducible if and only if it cannot be transformed into a block upper-triangular form by simultaneous row/column permutations.

\subsection{Problem Statements}

Consider a general model of opinion dynamics as follows:
\begin{equation}\label{sec2:1}
  \dot{x}_i=a_i(x_i)\bigg(-\sigma_ix_i+\sum_{j=1}^Nb_{ij}x_j\bigg),\;i=1,\cdots,N,
\end{equation}
where $x_i(t)\in[-1,1]$ is the agent state, $a_i(x_i)\in[0,1]$ is the susceptibility function which is state-dependent and Lipschitz in $[-1,1]$, and $\sigma_i>0$ represents forgetting factors for the $i$th agent. $B=[b_{ij}]\in\mathbb{R}^{N\times N}$ is the adjacency matrix of a strongly connected signed digraph $\mathcal{G}(B)$, where $b_{ii}=0$, $b_{ij}>0$ indicates that agents $i$ and $j$ are friendly, and $b_{ij}<0$ indicates that agents $i$ and $j$ are antagonistic.

Let $x=[x_1,\cdots,x_N]^{\rm T}$ and $\Sigma=\diag\{\sigma_1,\cdots,\sigma_N\}$. Then the generalized opinion dynamics (\ref{sec2:1}) can be written in matrix form as
\begin{equation}\label{sec2:2}
  \dot{x}=A(x)Ex,
\end{equation}
where $A(x)=\diag\{a_1(x_1),\cdots,a_N(x_N)\}$, and $E=B-\Sigma$.

If there do not exist the susceptibility functions, then the generalized opinion dynamics (\ref{sec2:2}) becomes $\dot{x}=Ex$, which was investigated in \cite{altafini2015predictable}. However, when there exist the susceptibility functions $a_i(x_i)\in[-1,1]$, the dynamic behavior of the generalized opinion dynamics (\ref{sec2:2}) will become more complex when compared to $\dot{x}=Ex$. If $\Sigma=I$ and the adjacency matrix $B\, (b_{ij}\in\{0,1\})$ is row-stochastic, then the generalized opinion dynamics (\ref{sec2:2}) becomes the model (1) considered in \cite{amelkin2017polar}. In \cite{amelkin2017polar}, the dynamic behavior of its model (1) was studied for some special susceptibility functions, that is, stubborn positives scenario $A(x)=0.5(I-\diag(x))$, stubborn neutrals scenario $A(x)=\diag(x)^2$, and stubborn extremists scenario $A(x)=(I-\diag(x)^2)$.

Our model (\ref{sec2:2}) can be viewed as a generalization of the models considered in \cite{altafini2015predictable,amelkin2017polar}. In the next section, we will study the dynamic behavior of the generalized opinion dynamics (\ref{sec2:2}) under three special scenarios, that is, stubborn extremists, stubborn positives, and stubborn neutrals. In this paper, we assume that matrix $E$ is eventually positive and all eigenvalues have non-positive real parts. Precisely, we adopt the following two assumptions.

\begin{assumption}\label{assumption1}
The adjacency matrix $B$ is eventually positive, and $\sigma_i=\sigma_j\geq\rho(B),\, \forall i,j=1,\cdots,N$.
\end{assumption}

\begin{assumption}\label{assumption2}
The adjacency matrix $B$ is not eventually positive, and $E$ can be represented as $E=C-dI$, where $C=B+D$ is eventually positive, $D=dI-\Sigma$, and $d\geq\rho(C)$.
\end{assumption}

\section{Analysis of Generalized Opinion Dynamics}\label{section3}

This section will study the dynamic behavior of the generalized opinion dynamics (\ref{sec2:2}) under three special scenarios (stubborn extremists, stubborn positives, and stubborn neutrals) respectively. Before moving forward, we first consider equilibrium points of the generalized opinion dynamics (\ref{sec2:2}).

All equilibrium points of the generalized opinion dynamics (\ref{sec2:2}) satisfy $A(x)Ex=0$. Suppose that $\mathcal{I}_{a,0}=\{i\mid a_i(x_i)=0\}$, $\mathcal{I}_{a,+}=\{i\mid a_i(x_i)>0\}$. First, if $\mathcal{I}_{a,0}=\varnothing$, then all $a_i(x_i)>0, i=1,\cdots,N$, and all equilibrium points satisfy $Ex=0$. Second, if $\mathcal{I}_{a,0}\neq\varnothing$ and $\mathcal{I}_{a,+}\neq\varnothing$, then there exists a permutation matrix $P$ such that
\begin{equation}\label{sec3:1}
  PA(x)P^{\rm T}=\left[\begin{array}{cc}A_1(x) & 0\\
  0 & A_2(x)\end{array}\right],
\end{equation}
where $A_1(x)=0, A_2(x)>0$. Let
\begin{equation*}
  y=Px=\left[\begin{array}{c}y_1\\
  y_2\end{array}\right],
\end{equation*}
where $\dim y_1=|\mathcal{I}_{a,0}|$. Then
\begin{align}\label{sec3:2}
  \dot{y}&=P\dot{x}=PA(x)Ex\notag\\
         &=PA(x)P^{\rm T}PEP^{\rm T}y\notag\\
         &=PA(x)P^{\rm T}\bar{E}y\notag\\
  &=\left[\begin{array}{cc}0 & 0\\
  0 & A_{2}(y_2)\end{array}\right]\left[\begin{array}{cc}\bar{E}_{11} & \bar{E}_{12}\\
  \bar{E}_{21} & \bar{E}_{22}\end{array}\right]\left[\begin{array}{c}y_1\\
  y_2\end{array}\right],
\end{align}
where $A_2(y_2)>0$. Hence, all equilibrium points of (\ref{sec3:2}) satisfy $y^{\ast}_1\in\{s\mid a_i(s)=0,\forall i\}^{|\mathcal{I}_{a,0}|}$, $\bar{E}_{22}y^{\ast}_2=-\bar{E}_{21}y^{\ast}_1$. In this case, all equilibrium points are $x^{\ast}=P^{\rm T}y^{\ast}$. Third, if $\mathcal{I}_{a,+}=\varnothing$, then all equilibrium points of the generalized opinion dynamics (\ref{sec2:2}) are $x^{\ast}\in\{s\mid a_i(s)=0,\forall i\}^N$.

\subsection{Stubborn Positives Scenario}\label{section31}

This subsection focuses on the dynamic behavior of the generalized opinion dynamics (\ref{sec2:2}) under the stubborn positives scenario, that is, $A(x)=0.5(I-\diag(x))$.

\begin{theorem}\label{theorem1}
Suppose that Assumption~\ref{assumption1} (Assumption~\ref{assumption2}) holds, $\mathcal{G}(B)$ is the signed digraph, and there exists a positive diagonal matrix $\Gamma=\diag\{\gamma_1,\cdots,\gamma_N\}$ such that $H(\Gamma E)\preceq0$ and $\rank H(\Gamma E)=\rank E$.
\begin{enumerate}[1)]
\item If $x(0)<\mathbf{1}$, then $\lim\limits_{t\rightarrow\infty}x(t)=\alpha v_r$, where $v_r$ is the right positive eigenvector of $B$ $(C)$, $\alpha\in\big[-\frac{1}{||v_r||_{\infty}},\frac{1}{||v_r||_{\infty}}\big)$.

\item If $\mathcal{I}_{a,0}\neq\varnothing$, $\mathcal{I}_{a,+}\neq\varnothing$, there exists a permutation matrix $P$ such that the network (\ref{sec2:2}) has the form of (\ref{sec3:2}), $\bar{E}_{22}$ is nonsingular, and $-\mathbf{1}\leq-\bar{E}_{22}^{-1}\bar{E}_{21}\mathbf{1}<\mathbf{1}$, then $\lim\limits_{t\rightarrow\infty}x(t)=x^{\ast}$ for all $x(0)\in \{x\mid x_i=1, i\in \mathcal{I}_{a,0}, x_i\in[-1,1), \forall i\notin \mathcal{I}_{a,0}\}$, where $x^{\ast}=P^{\rm T}y^{\ast}$, and $y^{\ast}$ is the equilibrium point of (\ref{sec3:2}).

\item If $\sigma_i=\sigma_j=\rho(B)$ $\left(d=\rho(C)\right)$, $\forall i,j=1,\cdots,N$, $v_r\neq\mathbf{1}$, then every solution starting in $\{x\mid x_i=1, \forall i\in \mathcal{J}, x_i\in[-1,1), \forall i\notin \mathcal{J}\}$ approaches $\beta v_r$, where $\mathcal{J}=\{j\mid v_{rj}=||v_r||_{\infty}\}$, $\beta=\frac{1}{||v_r||_{\infty}}$.
\end{enumerate}
\end{theorem}

\begin{IEEEproof}\
First, we consider the case that Assumption~\ref{assumption1} holds. The partition of the state space is illustrated in Fig.~\ref{figure1}.

1) Because $B$ is an eventually positive matrix, by Lemma~\ref{lemma2}, we have that matrix $B$ is irreducible, $\rho(B)$ is a simple positive eigenvalue of $B$ and the corresponding right eigenvector $v_r$ is positive. When Assumption~\ref{assumption1} holds, $\sigma_i=\sigma_j\geq\rho(B), \forall i,j=1,\cdots,N$. If $\sigma_i=\rho(B), \forall i=1,\cdots,N$, then all eigenvalues of $E=B-\Sigma$ have non-positive real parts. Moreover, $0$ is a simple eigenvalue of $E$ with a right positive eigenvector $v_r$. If $\sigma_i=\sigma_j>\rho(B), \forall i,j=1,\cdots,N$, then all eigenvalues of $E$ have negative real parts, and its dominant eigenvalue has the right positive eigenvector $v_r$.

\begin{figure}[t!]
\vspace*{-3ex}
\centering
\includegraphics[width=0.48\textwidth]{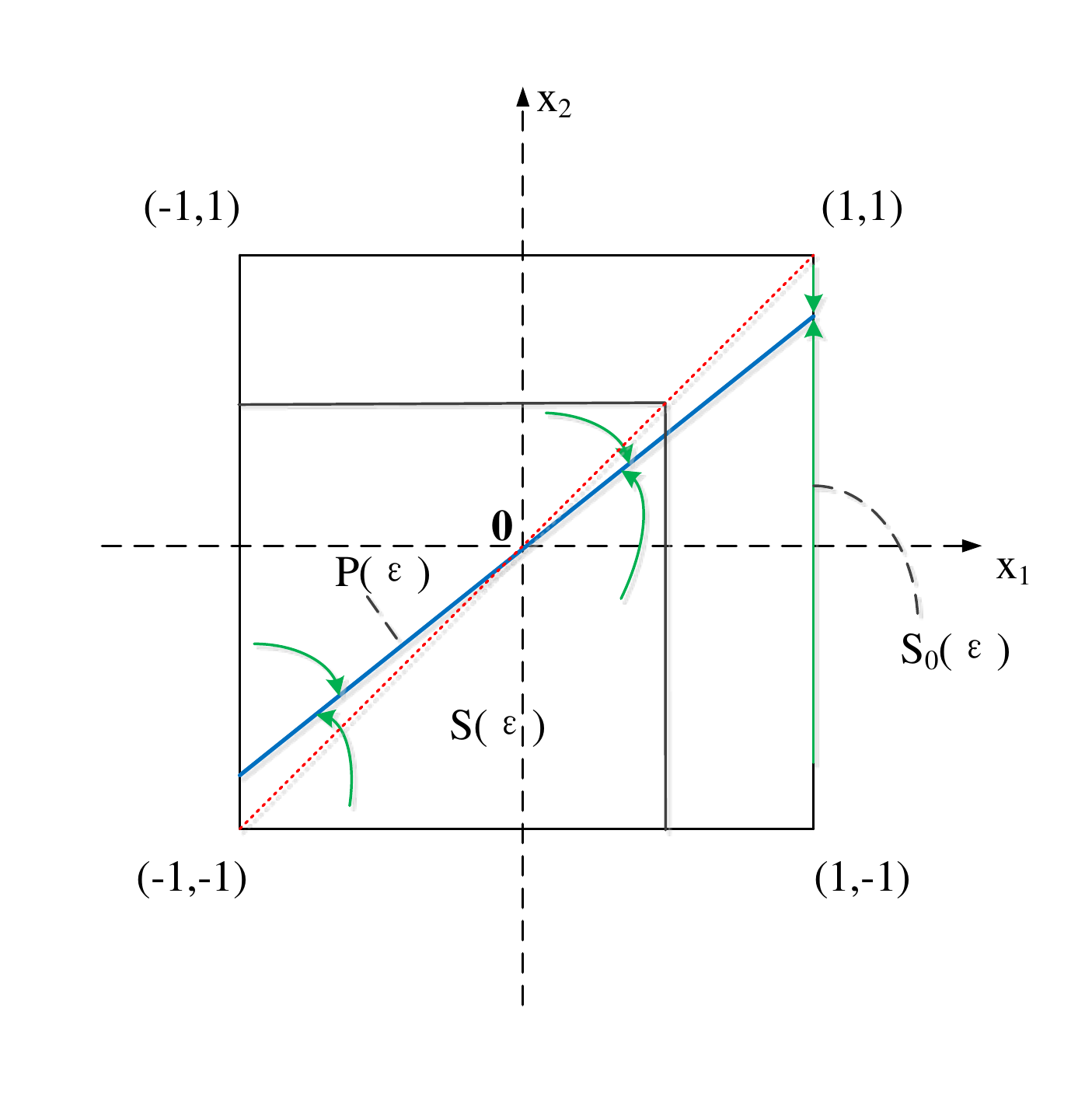}
\vspace{-4ex}
\caption{Convergence behavior of the generalized opinion dynamics (\ref{sec2:2}) with stubborn positives in two dimensions. The blue solid line is spanned by the right positive eigenvector $v_r$ and the red dotted line is spanned by vector $\mathbf{1}$.}
\label{figure1}
\end{figure}

Let $\mathcal{S}(\varepsilon)=[-1,1-\varepsilon]^N$, where $0<\varepsilon<2$. Suppose that
\begin{equation}\label{sec31:1}
  V(x)=\sum_{i=1}^N\gamma_i\left[-2x_i-2\ln(1-x_i)\right].
\end{equation}
Then $V(x)$ is positive definite in $\mathcal{S}(\varepsilon)$. For all $x\in \mathcal{S}(\varepsilon)$, the time derivative of $V(x)$ along (\ref{sec2:2}) is given by
\begingroup
\allowdisplaybreaks
\begin{align*}
  \dot{V}&=\sum_{i=1}^N2\gamma_ix_i\frac{\dot{x}_i}{1-x_i}\\
         &=\sum_{i=1}^Nx_i\gamma_ie_ix\\
         &=x^{\rm T}H(\Gamma E)x,
\end{align*}
\endgroup
where $e_i$ is the $i$th row of matrix $E$. Since $H(\Gamma E)\preceq0$, we have
\begin{align*}
  \dot{V}&=\sum_{i=1}^Nx_i\gamma_ie_ix\\
         &=x^{\rm T}H(\Gamma E)x\\
         &\leq0
\end{align*}
for all $x\in \mathcal{S}(\varepsilon)$. Hence, $\mathcal{S}(\varepsilon)$ is positively invariant with respect to (\ref{sec2:2}). By Lemma~\ref{lemma4}, $\dot{V}=0$ if and only if $H(\Gamma E)x=\mathbf{0}$. Because $\rank H(\Gamma E)=\rank E$, $H(\Gamma E)$ and $E$ have the same null space. Thus, $\dot{V}=0$ if and only if $Ex=\mathbf{0}$. If $\sigma_i=\rho(B), \forall i=1,\cdots,N$, then $\dot{V}=0$ if and only if $x\in \mathcal{P}(\varepsilon)=\big\{x\;\big|\; x=\alpha v_r, \alpha\in\big[-\frac{1}{||v_r||_{\infty}},\frac{1-\varepsilon}{||v_r||_{\infty}}\big]\big\}$. If $\sigma_i=\sigma_j>\rho(B), \forall i,j=1,\cdots,N$, then $\dot{V}=0$ if and only if $x=\mathbf{0}$. By LaSalle invariance principle, the solution of (\ref{sec2:2}) with initial condition $x(0)\in \mathcal{S}(\varepsilon)$ converges into $\mathcal{P}(\varepsilon)$. If $x(0)$ is sufficiently close to $\mathbf{1}$, then $\varepsilon\rightarrow 0$. Therefore, for all $x(0)<\mathbf{1}$, $\lim\limits_{t\rightarrow\infty}x(t)=\alpha v_r,\, \alpha\in\big[-\frac{1}{||v_r||_{\infty}},\frac{1}{||v_r||_{\infty}}\big)$.

2) Let $\mathcal{S}_{\mathcal{I}_{a,0}}(\varepsilon)=\{x\mid x_i=1, \forall i\in \mathcal{I}_{a,0},\; x_i\in[-1,1-\varepsilon], \forall i\notin \mathcal{I}_{a,0}\}$, where $0<\varepsilon<2$. When $x(0)\in \mathcal{S}_{\mathcal{I}_{a,0}}(\varepsilon)$, the solution of (\ref{sec2:2}) satisfies $x_i(t)\equiv 1, \forall i\in \mathcal{I}_{a,0}$. $P$ is a permutation matrix such that
\begin{equation*}
  y=Px=\left[\begin{array}{c}y_1\\
  y_2\end{array}\right],
\end{equation*}
where $y_1=\mathbf{1}, \dim y_1=|\mathcal{I}_{a,0}|$. Then the network (\ref{sec2:2}) becomes
\begin{align*}
  \dot{y}
  &=\left[\begin{array}{cc}0 & 0\\
  0 & 0.5(I-\diag(y_2))\end{array}\right]\left[\begin{array}{cc}\bar{E}_{11} & \bar{E}_{12}\\
  \bar{E}_{21} & \bar{E}_{22}\end{array}\right]\left[\begin{array}{c}y_1\\
  y_2\end{array}\right].
\end{align*}
Hence,
\begin{equation}\label{sec31:2}
  \dot{y}_2=0.5(I-\diag(y_2))\left(\bar{E}_{22}y_2+\bar{E}_{21}\mathbf{1}\right),
\end{equation}
where $y_2\in[-1,1-\varepsilon]^{|\mathcal{I}_{a,+}|}$.

Since $\bar{E}_{22}$ is nonsingular, we have that $\det(\bar{E}_{22})\neq 0$ and there exists a unique $y_2^{\ast}$ such that $\bar{E}_{22}y_2^{\ast}=-\bar{E}_{21}\mathbf{1}$. Let $s=y_2-y_2^{\ast}$. Then we obtain
\begin{equation}\label{sec31:3}
  \dot{s}=0.5(I-\diag(s+y_2^{\ast}))\bar{E}_{22}s.
\end{equation}
where $s\in[-1-y_2^{\ast},1-\varepsilon-y_2^{\ast}]^{|\mathcal{I}_{a,+}|}$. Assume that
\begin{align}\label{sec31:4}
  V(s)
  =-\sum_{i=1}^{\bar{N}}\bar{\gamma}_{i+|\mathcal{I}_{a,0}|}\left[2s_i+2(1-y_{2i}^{\ast})\ln(1-y_{2i}^{\ast}-s_i)\right],
\end{align}
where $\bar{N}=N-|\mathcal{I}_{a,0}|$, $\bar{\Gamma}=\diag\{\bar{\gamma}_1,\cdots,\bar{\gamma}_N\}=P\Gamma P^{\rm T}$. Let $\bar{\Gamma}_1=\diag\{\bar{\gamma}_1,\cdots,\bar{\gamma}_{|\mathcal{I}_{a,0}|}\}$ and $\bar{\Gamma}_2=\diag\{\bar{\gamma}_{|\mathcal{I}_{a,0}|+1},\cdots,\bar{\gamma}_N\}$. It is easy to see from (\ref{sec31:4}) that $V(s)$ is continuously differentiable in  $[-1-y_2^{\ast},1-\varepsilon-y_2^{\ast}]^{|\mathcal{I}_{a,+}|}$. The time derivative of $V(s)$ along (\ref{sec31:3}) is calculated as follows:
\begin{align*}
  \dot{V}&=\sum_{i=1}^{N-|\mathcal{I}_{a,0}|}2\bar{\gamma}_{i+|\mathcal{I}_{a,0}|}\frac{s_i\dot{s}_i}{1-y_{2i}^{\ast}-s_i}\\
         &=\sum_{i=1}^{N-|\mathcal{I}_{a,0}|}s_i\bar{\gamma}_{i+|\mathcal{I}_{a,0}|}\bar{e}_{22i}s\\
         &=s^{\rm T}H(\bar{\Gamma}_2\bar{E}_{22})s,
\end{align*}
where $\bar{e}_{22i}$ is the $i$th row of matrix $\bar{E}_{22}$. Since $H(\Gamma E)\preceq0$ and $H(\bar{\Gamma}_2\bar{E}_{22})$ is a principal submatrix of $H(\Gamma E)$, we have $H(\bar{\Gamma}_2\bar{E}_{22})\preceq0$.
Hence,
\begin{align*}
  \dot{V}&=\sum_{i=1}^{N-|\mathcal{I}_{a,0}|}s_i\bar{\gamma}_{i+|\mathcal{I}_{a,0}|}\bar{e}_{22i}s\\
         &=s^{\rm T}H(\bar{\Gamma}_2\bar{E}_{22})s\\
         &\leq 0
\end{align*}
for all $s\in [-1-y_2^{\ast},1-\varepsilon-y_2^{\ast}]^{|\mathcal{I}_{a,+}|}$.

Thus, $[-1-y_2^{\ast},1-\varepsilon-y_2^{\ast}]^{|\mathcal{I}_{a,+}|}$ is positively invariant with respect to (\ref{sec31:3}). By LaSalle invariance principle, the solution of (\ref{sec31:3}) with initial condition $s(0)\in [-1-y_2^{\ast},1-\varepsilon-y_2^{\ast}]^{|\mathcal{I}_{a,+}|}$ converges to $y_2^{\ast}$. When $\varepsilon\rightarrow 0^{+}$, we have $1-\varepsilon\rightarrow1^{-}$. Hence, if $x(0)\in \{x\mid x_i=1, i\in \mathcal{I}_{a,0},\, x_i\in[-1,1), \forall i\notin \mathcal{I}_{a,0}\}$, then $\lim\limits_{t\rightarrow\infty}x(t)=x^{\ast}$.

3) Let $\mathcal{S}_{0}(\varepsilon)=\{x\mid x_i=1, \forall i\in \mathcal{J}, x_i\in[-1,1-\varepsilon], \forall i\notin \mathcal{J}\}$, where $0<\varepsilon<2$, $\mathcal{J}=\{j\mid v_{rj}=||v_r||_{\infty}\}$. There exists a permutation matrix $P$ such that the network (\ref{sec2:2}) takes the following form
\begin{align*}
  \dot{y}
  &=\left[\begin{array}{cc}0 & 0\\
  0 & 0.5(I-\diag(y_2))\end{array}\right]\left[\begin{array}{cc}\bar{E}_{11} & \bar{E}_{12}\\
  \bar{E}_{21} & \bar{E}_{22}\end{array}\right]\left[\begin{array}{c}y_1\\
  y_2\end{array}\right],
\end{align*}
where
\begin{equation*}
  y=Px=\left[\begin{array}{c}y_1\\
  y_2\end{array}\right]
\end{equation*}
with $y_1=\mathbf{1}, \dim y_1=|\mathcal{J}|$. Then we can get
\begin{equation}\label{sec31:5}
  \dot{y}_2=0.5(I-\diag(y_2))\left(\bar{E}_{22}y_2+\bar{E}_{21}\mathbf{1}\right).
\end{equation}
Because $\sigma_i=\sigma_j=\rho(B)$ ($d=\rho(C)$), $\forall i,j=1,\cdots,N$, $v_r\neq\mathbf{1}$ and $\mathcal{J}=\{j\mid v_{rj}=||v_r||_{\infty}\}$, there exists a unique solution $y_2^{\ast}$ such that $\bar{E}_{22}y_2^{\ast}+\bar{E}_{21}\mathbf{1}=0$. From a geometric point of view, the point
\begin{equation*}
  x^{\ast}=P^{\rm T}\left[\begin{array}{c}\mathbf{1}\\
  y_2^{\ast}\end{array}\right]
\end{equation*}
is the intersection point of the line $\alpha v_r$ with $[-1,1]^N$. Similar to the arguments presented in 2), we can obtain that every solution starting in $\mathcal{S}_{0}(\varepsilon)$ approaches $x^{\ast}= \beta v_r$, where $\beta=\frac{1}{||v_r||_{\infty}}$. When $\varepsilon\rightarrow 0^{+}$, we can get that every solution starting in $\{x\mid x_i=1, \forall i\in \mathcal{J},\, x_i\in[-1,1), \forall i\notin \mathcal{J}\}$ tends to $\beta v_r$.

Second, we consider the case that Assumption~\ref{assumption2} holds. Because $C=B+D$ is eventually positive, by  Lemma~\ref{lemma2}, matrix $C$ is irreducible, $\rho(C)$ is a simple positive eigenvalue of $C$ and the corresponding right eigenvector $v_r$ is positive. Hence, $E=C-dI$ has the eigenvalue $\rho(C)-d$ with an eigenvector $v_r$. If $d=\rho(C)$, then all eigenvalues of $E$ have non-positive real parts, and $0$ is a simple eigenvalue of $E$ with a right positive eigenvector $v_r$. If $d>\rho(C)$, then all eigenvalues of $E$ have negative real parts, and its dominant eigenvalue has the right positive eigenvector $v_r$.

On this basis, one can use the similar methods presented above to acquire the results in 1), 2) and 3) in the case of Assumption~\ref{assumption2}, and details are omitted here for brevity. This completes the proof of Theorem~\ref{theorem1}.
\end{IEEEproof}

\begin{corollary}\label{corollary1}
Suppose that Assumption~\ref{assumption1} (Assumption~\ref{assumption2}) holds and the signed graph $\mathcal{G}(B)$ is undirected.
\begin{enumerate}[1)]
\item If $x(0)<\mathbf{1}$, then $\lim\limits_{t\rightarrow\infty}x(t)=\alpha v_r$, where $v_r$ is the right positive eigenvector of $B$ $(C)$, $\alpha\in\big[-\frac{1}{||v_r||_{\infty}},\frac{1}{||v_r||_{\infty}}\big)$.

\item If $\mathcal{I}_{a,0}\neq\varnothing$, $\mathcal{I}_{a,+}\neq\varnothing$, there exists a permutation matrix $P$ such that the network (\ref{sec2:2}) has the form of (\ref{sec3:2}), $\bar{E}_{22}$ is nonsingular, and $-\mathbf{1}\leq-\bar{E}_{22}^{-1}\bar{E}_{21}\mathbf{1}<\mathbf{1}$, then $\lim\limits_{t\rightarrow\infty}x(t)=x^{\ast}$ for all $x(0)\in \{x\mid x_i=1, i\in \mathcal{I}_{a,0},\, x_i\in[-1,1), \forall i\notin \mathcal{I}_{a,0}\}$, where $x^{\ast}=P^{\rm T}y^{\ast}$, and $y^{\ast}$ is the equilibrium points of (\ref{sec3:2}).

\item If $\sigma_i=\sigma_j=\rho(B)$ $\left(d=\rho(C)\right)$, $\forall i,j=1,\cdots,N$, $v_r\neq\mathbf{1}$, then every solution starting in $\{x\mid x_i=1, \forall i\in \mathcal{J}, x_i\in[-1,1), \forall i\notin \mathcal{J}\}$ approaches $\beta v_r$, where $\mathcal{J}=\{j\mid v_{rj}=||v_r||_{\infty}\}$, $\beta=\frac{1}{||v_r||_{\infty}}$.
\end{enumerate}
\end{corollary}

\begin{IEEEproof}\
Because the signed graph $\mathcal{G}(B)$ is undirected, it can be known that matrix $E$ is symmetric and diagonalizable. By Remark~\ref{remark3} given in Subsection~\ref{section32}, we have $H(E)\preceq0$ and $\rank H(E)=\rank E$. Hence, it follows from Theorem~\ref{theorem1} that 1), 2) and 3) of this corollary hold. This completes the proof of Corollary~\ref{corollary1}.
\end{IEEEproof}

The stubborn positives scenario describes the case that the agents only at one end of the opinion spectrum are stubborn. For this scenario, an example about two smartphone brands was presented in \cite{amelkin2017polar}. The opinion values $1$ and $-1$ mean an aggressive marketing of the brand and a neutral marketing of the brand, respectively. Opinion $-1$ may be changed by other opinions, but opinion $1$ is stubborn and cannot be influenced by other opinions. The above Theorem~\ref{theorem1} and Corollary~\ref{corollary1} show that the states of the generalized opinion dynamics (\ref{sec2:2}) are convergent under different initial conditions.

\subsection{Stubborn Neutrals Scenario}\label{section32}

This subsection examines the dynamic behavior of the generalized opinion dynamics (\ref{sec2:2}) under the stubborn neutrals scenario, that is, $A(x)=\diag(x)^2$.

\begin{lemma}\label{lemma1}
Consider the generalized opinion dynamics (\ref{sec2:2}) under the stubborn neutrals scenario $A(x)=\diag(x)^2$. Then $x(t)$ has nonnegative derivatives.
\end{lemma}

\begin{IEEEproof}\
Let $x_t(z)$ be the solution of the generalized opinion dynamics (\ref{sec2:2}) with initial condition $z\in\mathbb{R}^N$. Define $M(t)=Dx_t(z)$, $F(x)=\diag(x)^2Ex$, where $Dx_t(z)$ is the Jacobian matrix. Then $M(t)$ satisfies the variational equation
\begin{equation*}
  \dot{M}=DF(x)M,\;M(0)=I,
\end{equation*}
where $DF(x)$ is the Jacobian matrix. It is easy to get
\begin{equation*}
  DF(x)=\diag(x)\left[\begin{array}{cccc}
         \eta _1 & e_{12}x_1 & \cdots & e_{1N}x_1\\
        e_{21}x_2 & \eta_2 & \cdots & e_{2N}x_2\\
        \cdots & \cdots & \cdots & \cdots\\
        e_{N1}x_N & e_{N2}x_N & \cdots & \eta_N\\
        \end{array}\right],
\end{equation*}
where $\eta_i=e_{ii}x_i+2\sum_{j=1}^Ne_{ij}x_j$, $i=1,2,\cdots,N$. We can obtain that $M(t)=0$ when $x=\mathbf{0}$. Hence, $M(t)\geq0$, that is, $x(t)$ has nonnegative derivatives.
\end{IEEEproof}

\begin{theorem}\label{theorem2}
Suppose that Assumption~\ref{assumption1} (Assumption~\ref{assumption2}) holds, $\mathcal{G}(B)$ is the signed digraph, and there exists a positive diagonal matrix $\Gamma=\diag\{\gamma_1,\cdots,\gamma_N\}$ such that $H(\Gamma E)\preceq0$ and $\rank H(\Gamma E)=\rank E$.
If $\sigma_i=\sigma_j=\rho(B)$ $\left(d=\rho(C)\right)$, $\forall i,j=1,\cdots,N$, then
\begin{enumerate}[1)]
\item $\lim\limits_{t\rightarrow\infty}x(t)=\alpha v_r$ when $x(0)>\mathbf{0}$, where $v_r$ is the right positive eigenvector of $B$ $(C)$, $\alpha\in\big(0,\frac{1}{||v_r||_{\infty}}\big]$;

\item $\lim\limits_{t\rightarrow\infty}x(t)=\alpha v_r$ when $x(0)<\mathbf{0}$, where $\alpha\in\big[-\frac{1}{||v_r||_{\infty}},0\big)$;

\item $\lim\limits_{t\rightarrow\infty}x(t)=\mathbf{0}$ otherwise.
\end{enumerate}
If $\sigma_i=\sigma_j>\rho(B)$ $\left(d>\rho(C)\right)$, $\forall i,j=1,\cdots,N$, then $\lim\limits_{t\rightarrow\infty}x(t)=\mathbf{0},\, \forall x(0)\in [-1,1]^N$.
\end{theorem}

\begin{IEEEproof}\
First, we consider the case that Assumption~\ref{assumption1} holds. The corresponding state space is partitioned as depicted in Fig.~\ref{figure2}.

\begin{figure}[t!]
\vspace*{-3ex}
\centering
\includegraphics[width=0.48\textwidth]{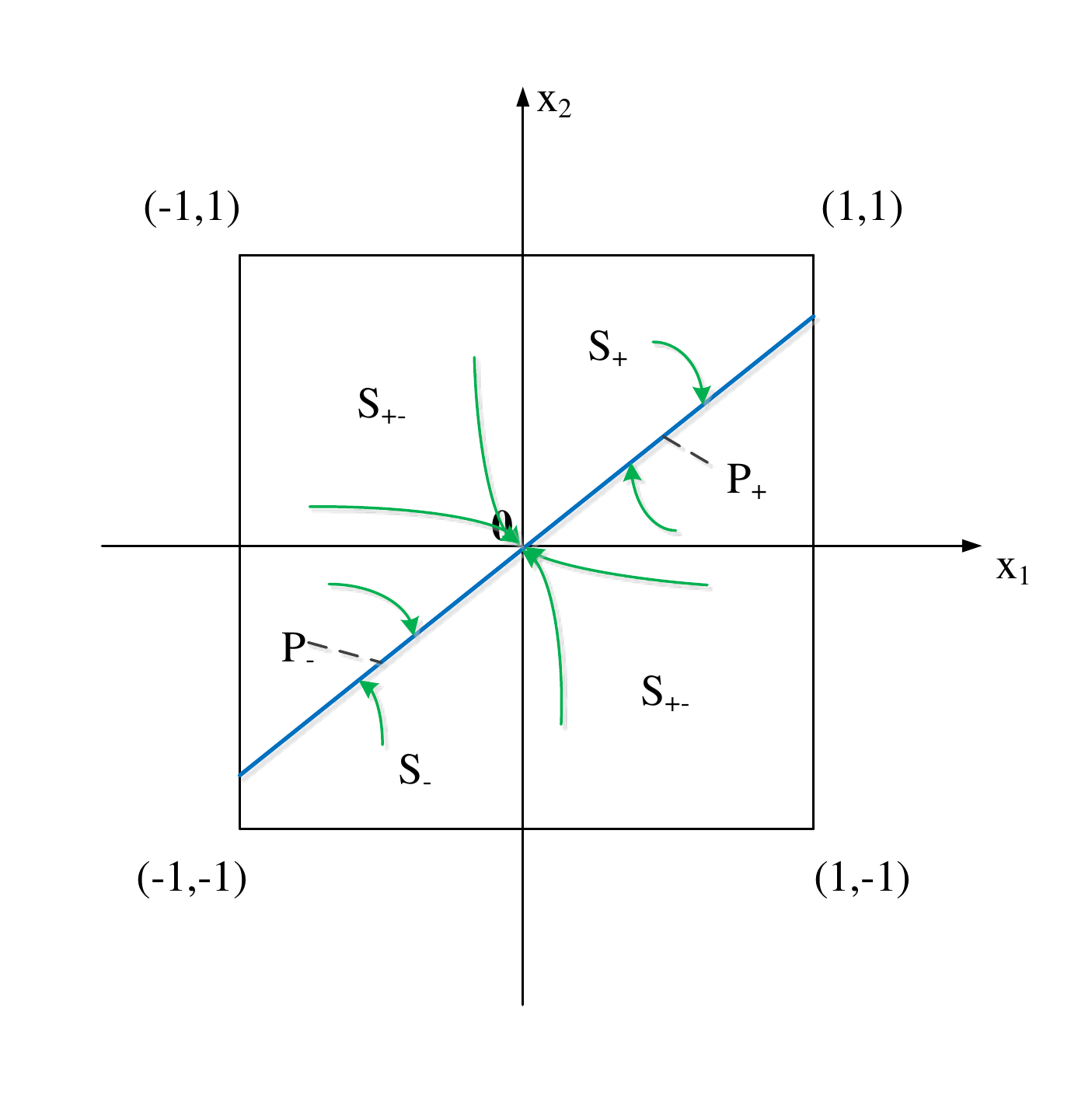}
\vspace{-4ex}
\caption{Convergence behavior of the generalized opinion dynamics (\ref{sec2:2}) with stubborn neutrals in two dimensions. The blue solid line is spanned by the right positive eigenvector $v_r$, and the green curves with arrows denote several trajectories of the generalized opinion dynamics (\ref{sec2:2}).}
\label{figure2}
\end{figure}

When $\sigma_i=\sigma_j=\rho(B),\, \forall i,j=1,\cdots,N$, and $B$ is eventually positive matrix, by Lemma~\ref{lemma2}, we have that matrix $B$ is irreducible, $\rho(B)$ is a simple positive eigenvalue of $B$ and the corresponding right eigenvector $v_r$ is positive. In this case, the equilibrium points are $x^{\ast}=\alpha v_r$, $\alpha\in\big[-\frac{1}{||v_r||_{\infty}},\frac{1}{||v_r||_{\infty}}\big]$.
Specifically, if $\mathcal{I}_{a,0}=\varnothing$, then all equilibrium points satisfy $Ex=0$. Hence, the equilibrium points are $x^{\ast}=\alpha v_r$. If $\mathcal{I}_{a,+}=\varnothing$, then the equilibrium points are $x^{\ast}=\mathbf{0}$. If $\mathcal{I}_{a,0}\neq\varnothing$ and $\mathcal{I}_{a,+}\neq\varnothing$, then there exists a permutation matrix $P$ such that
\begin{align*}
  \bar{E}&=PEP^{\rm T},\\
  &=\left[\begin{array}{cc}\bar{E}_{11} & \bar{E}_{12}\\
  \bar{E}_{21} & \bar{E}_{22}\end{array}\right],
\end{align*}
and all equilibrium points are $x^{\ast}=P^{\rm T}y^{\ast}$, $y^{\ast}_1=\mathbf{0}$, $\bar{E}_{22}y^{\ast}_2=\mathbf{0}$. If $\bar{E}_{22}$ is nonsingular, then $y^{\ast}_2=\mathbf{0}$ and $x^{\ast}=\mathbf{0}$. Assume that $\bar{E}_{22}$ is singular, and there exists a nonzero vector $\xi$ such that $\bar{E}_{22}\xi=\mathbf{0}$. Let $x=[\mathbf{0}^{\rm T}\; \xi^{\rm T}]^{\rm T}$. Since there exists a positive diagonal matrix $\Gamma$ such that $H(\Gamma E)\preceq0$, and $\rank H(\Gamma E)=\rank E$, $H(\Gamma E)$ and $E$ have the same null space. As $P$ is a permutation matrix, $H(P\Gamma EP^{\rm T})$ and $PEP^{\rm T}$ have the same null space. Suppose that $\bar{\Gamma}=\diag\{\bar{\Gamma}_1,\bar{\Gamma}_2\}=P\Gamma P^{\rm T}$. Then $x^{\rm T}H(P\Gamma EP^{\rm T})x=\xi^{\rm T}H(\bar{\Gamma}_2\bar{E}_{22})\xi=0$. Thus, $PEP^{\rm T}x=\mathbf{0}$. Because $x=[\mathbf{0}^{\rm T}\; \xi^{\rm T}]^{\rm T}$, $PEP^{\rm T}x=\mathbf{0}$ if and only if $\xi=\mathbf{0}$. This is a contradiction with the assumption that $\xi$ is a nonzero vector. Therefore, if $\mathcal{I}_{a,0}\neq\varnothing$ and $\mathcal{I}_{a,+}\neq\varnothing$, then the equilibrium points are $x^{\ast}=\mathbf{0}$.

If $\sigma_i=\sigma_j>\rho(B), \forall i,j=1,\cdots,N$, then we can use the same method to obtain that the equilibrium points are $x^{\ast}=\mathbf{0}$.

1) If $x(0)>\mathbf{0}$, then there exists $\varepsilon\in(0,1)$ such that $x(0)\in[\varepsilon,1]^N$. Let
\begin{equation}\label{sec32:1}
  V(x)=2\sum_{i=1}^N\gamma_i\ln(1+x_i).
\end{equation}
So $V(x)$ is continuously differentiable in $\mathcal{S}_{+}=[\varepsilon,1]^N$. By Lemma~\ref{lemma1}, $x(t)$ has nonnegative derivatives, that is, $\dot{x}_i\geq0, \forall i=1,\cdots,N$. For all $x\in \mathcal{S}_{+}$, the time derivative of $V(x)$ along (\ref{sec2:2}) is given by
\begin{align*}
  \dot{V}&=2\sum_{i=1}^N\gamma_i\frac{\dot{x}_i}{1+x_i}\\
         &\leq 2\sum_{i=1}^N\gamma_i\frac{\dot{x}_i}{2x_i}\\
         &=x^{\rm T}H(\Gamma E)x\\
         &\leq0
\end{align*}
for all $x\in \mathcal{S}_{+}$. Hence, $\mathcal{S}_{+}$ is positively invariant with respect to (\ref{sec2:2}). $\dot{V}=0$ if and only if $x\in \mathcal{P}_{+}=\big\{x \,\big|\, x=\alpha v_r, \alpha\in\big[\frac{\varepsilon}{||v_r||_{\infty}}, \frac{1}{||v_r||_{\infty}}\big]\big\}$. By LaSalle invariance principle, the solution of (\ref{sec2:2}) with initial condition $x(0)\in \mathcal{S}_{+}$ converges into $\mathcal{P}_{+}$. Hence, if $x(0)>\mathbf{0}$, then $\lim\limits_{t\rightarrow\infty}x(t)=\alpha v_r$, $\alpha\in\big(0,\frac{1}{||v_r||_{\infty}}\big]$.

2) If $x(0)<\mathbf{0}$, then there exists $\varepsilon\in(0,1)$ such that $x(0)\in[-1,-\varepsilon]^N$. Let
\begin{equation}\label{sec32:2}
  V(x)=2\sum_{i=1}^N\gamma_i\ln(\varepsilon-x_i).
\end{equation}
So $V(x)$ is continuously differentiable in $\mathcal{S}_{-}=[-1,-\varepsilon]^N$. Because $x(t)$ has nonnegative derivatives according to Lemma~\ref{lemma1}, the time derivative of $V(x)$ along (\ref{sec2:2}) is found to be
\begin{align*}
  \dot{V}&=2\sum_{i=1}^N\gamma_i\frac{\dot{x}_i}{x_i-\varepsilon}\\
         &\leq2\sum_{i=1}^N\gamma_i\frac{\dot{x}_i}{2x_i}\\
         &=x^{\rm T}H(\Gamma E)x\\
         &\leq0
\end{align*}
for all $x\in \mathcal{S}_{-}$. Hence, $\mathcal{S}_{-}$ is positively invariant with respect to (\ref{sec2:2}). $\dot{V}=0$ if and only if $x\in \mathcal{P}_{-}=\big\{x \,\big|\, x=\alpha v_r, \alpha\in\big[-\frac{1}{||v_r||_{\infty}},-\frac{\varepsilon}{||v_r||_{\infty}}\big]\big\}$. From LaSalle invariance principle, it follows that the solution of (\ref{sec2:2}) with initial condition $x(0)<\mathbf{0}$ converges into $\mathcal{P}_{-}$.

3) Let $\mathcal{S}_{\pm}=\big\{x \,\big|\, x\in[-1,1]^N, \prod\limits_{i=1}^Nx_i=0\big\}\bigcup\big\{x \,\big|\, x\in[-1,1]^N, \exists i,j \;\text{s.t.}\; \sgn(x_ix_j)=-1\big\}$. By Lemma~\ref{lemma1}, $x(t)$ has nonnegative derivatives, that is, $\dot{x}_i\geq0, \forall i=1,\cdots,N$. From the generalized opinion dynamics (\ref{sec2:2}) with $A(x)=\diag(x)^2$, it can be seen that $\sgn(x_i(t))=\sgn(x_i(0))$. For arbitrary $x(0)\in \mathcal{S}_{\pm}$ with $x(0)\neq \mathbf{0}$, without loss of generality, assume that $x_i(0)>0, i=1,2,\cdots,l_1$, $x_i(0)=0, i=l_1+1,\cdots,l_2$, $x_i(0)<0, i=l_2+1,\cdots,N$, where $1<l_1<l_2<N$ are positive constants. Define
\begin{equation}\label{sec32:3}
  V(x)=2\sum_{i=1}^{l_2}\gamma_i\ln(x_i+1)+2\sum_{i=l_2+1}^N\gamma_i\ln(\varepsilon-x_i),
\end{equation}
where $x(t)$ is the solution of the generalized opinion dynamics (\ref{sec2:2}) with initial condition $x(0)$, and $\varepsilon$ is a sufficiently small positive constant. Let $\mathcal{S}(\varepsilon)=\{x\mid x\in[-1,1]^N, x_i\in[-1,-\varepsilon], i=l_2+1,\cdots,N\}$. Then $V(x)$ is continuously differentiable in $\mathcal{S}_{\pm}\bigcap \mathcal{S}(\varepsilon)$. The time derivative of $V(x)$ along (\ref{sec2:2}) is given by
\begin{align*}
  \dot{V}&=2\sum_{i=1}^{l_1}\gamma_i\frac{\dot{x}_i}{1+x_i}+2\sum_{i=l_2+1}^N\gamma_i\frac{\dot{x}_i}{x_i-\varepsilon}\\
         &\leq 2\sum_{i=1}^{l_1}\gamma_i\frac{\dot{x}_i}{2x_i}+2\sum_{i=1}^N\gamma_i\frac{\dot{x}_i}{2x_i}\\
         &=x^{\rm T}H(\Gamma E)x
\end{align*}
for all $x\in \mathcal{S}_{\pm}\bigcap \mathcal{S}(\varepsilon)$. Since $x^{\rm T}H(\Gamma E)x=0, x\in \mathcal{S}_{\pm}$ if and only if $x=\mathbf{0}$, $\dot{V}<0$ for all $x\in \mathcal{S}_{\pm}\bigcap \mathcal{S}(\varepsilon)$. Hence, $V(x)$ is strictly monotonically decreasing on $\mathcal{S}_{\pm}\bigcap \mathcal{S}(\varepsilon)$, and $V(x)$ has a minimum value at $\bar{x}$ with $\bar{x}_i=0, i=1,\cdots,l_2, \bar{x}_i=-\varepsilon, i=l_2+1,\cdots,N$. Moreover, $\dot{V}(\bar{x})<0$. Thus, the solution of (\ref{sec2:2}) with initial condition $x(0)$ converges to $\mathbf{0}$. For arbitrary $x(0)\in \mathcal{S}_{\pm}$ , $\lim\limits_{t\rightarrow\infty}x(t)=\mathbf{0}$.

If $\sigma_i=\sigma_j>\rho(B), \forall i,j=1,\cdots,N$, then $H(E)\prec0$, and $Ex=\mathbf{0}$ if and only if $x=\mathbf{0}$. Hence, using the method shown above, it is easy to obtain that $\lim\limits_{t\rightarrow\infty}x(t)=\mathbf{0}, \forall x(0)\in [-1,1]^N$.

Second, we turn to the case that Assumption~\ref{assumption2} holds.
Since $E=C-dI$, $C=B+D$ is eventually positive, and $d\geq\rho(C)$, with the above-presented similar method we arrive at that the equilibrium points are $x^{\ast}=\alpha v_r$, $\alpha\in\big[-\frac{1}{||v_r||_{\infty}},\frac{1}{||v_r||_{\infty}}\big]$.
Precisely, if $\mathcal{I}_{a,0}=\varnothing$, then the equilibrium points are $x^{\ast}=\alpha v_r$. If $\mathcal{I}_{a,+}=\varnothing$, then the equilibrium points are $x^{\ast}=\mathbf{0}$. If $\mathcal{I}_{a,0}\neq\varnothing$ and $\mathcal{I}_{a,+}\neq\varnothing$, then all equilibrium points are $x^{\ast}=\mathbf{0}$. Therefore, the results of this theorem under Assumption~\ref{assumption2} can be derived by means of the method given above, but the actual derivation is skipped here. The proof of Theorem~\ref{theorem2} is thus complete.
\end{IEEEproof}

When the signed graph $\mathcal{G}(B)$ is undirected, we have that matrix $E$ is symmetric and diagonalizable. Then we can obtain the following corollary from Theorem~\ref{theorem2}. Its proof is similar to that of Corollary~\ref{corollary1}, so it is omitted here.

\begin{corollary}\label{corollary2}
Suppose that Assumption~\ref{assumption1} (Assumption~\ref{assumption2}) holds and the signed graph $\mathcal{G}(B)$ is undirected.
If $\sigma_i=\sigma_j=\rho(B)$ $\left(d=\rho(C)\right)$, $\forall i,j=1,\cdots,N$, then
\begin{enumerate}[1)]
\item $\lim\limits_{t\rightarrow\infty}x(t)=\alpha v_r$ when $x(0)>\mathbf{0}$, where $v_r$ is the right positive eigenvector of $B$ $(C)$, $\alpha\in\big(0,\frac{1}{||v_r||_{\infty}}\big]$;

\item $\lim\limits_{t\rightarrow\infty}x(t)=\alpha v_r$ when $x(0)<\mathbf{0}$, where $\alpha\in\big[-\frac{1}{||v_r||_{\infty}},0\big)$;

\item $\lim\limits_{t\rightarrow\infty}x(t)=\mathbf{0}$ otherwise.
\end{enumerate}
If $\sigma_i=\sigma_j>\rho(B)$ $\left(d>\rho(C)\right)$, $\forall i,j=1,\cdots,N$, then $\lim\limits_{t\rightarrow\infty}x(t)=\mathbf{0},\, \forall x(0)\in [-1,1]^N$.
\end{corollary}

The stubborn neutrals scenario represents the situation that the neutral opinions cannot be altered, while the extreme opinions may be easily affected by other opinions. Theorem~\ref{theorem2} and Corollary~\ref{corollary2} establish the convergence of the generalized opinion dynamics (\ref{sec2:2}) under different initial conditions for this scenario.

\subsection{Stubborn Extremists Scenario}\label{section33}

This subsection concentrates on the dynamic behavior of the generalized opinion dynamics (\ref{sec2:2}) under the stubborn extremists scenario, that is, $A(x)=(I-\diag(x)^2)$.

\begin{theorem}\label{theorem3}
Suppose that Assumption~\ref{assumption1} (Assumption~\ref{assumption2}) holds, $\mathcal{G}(B)$ is the signed digraph, and there exists a positive diagonal matrix $\Gamma=\diag\{\gamma_1,\cdots,\gamma_N\}$ such that $H(\Gamma E)\preceq0$ and $\rank H(\Gamma E)=\rank E$.
\begin{enumerate}[1)]
\item If $x(0)\in (-1,1)^N$, then $\lim\limits_{t\rightarrow\infty}x(t)=\alpha v_r$, where $v_r$ is the right positive eigenvector of $B$ $(C)$, and $\alpha\in\big(-\frac{1}{||v_r||_{\infty}},\frac{1}{||v_r||_{\infty}}\big)$.

\item If $\mathcal{I}_{a,0}\neq\varnothing$, $\mathcal{I}_{a,+}\neq\varnothing$, there exists a permutation matrix $P$ such that the network (\ref{sec2:2}) has the form of (\ref{sec3:2}), $\bar{E}_{22}$ is nonsingular, and $-\mathbf{1}<-\bar{E}_{22}^{-1}\bar{E}_{21}y_1^{\ast}<\mathbf{1}$, then $\lim\limits_{t\rightarrow\infty}x(t)=x^{\ast}$ for all $x(0)\in \{x\mid |x_i|=1, i\in \mathcal{I}_{a,0},\, x_i\in(-1,1), \forall i\notin \mathcal{I}_{a,0}\}$, where $x^{\ast}=P^{\rm T}y^{\ast}$, and $y^{\ast}$ is the equilibrium point of (\ref{sec3:2}).

\item If $\sigma_i=\sigma_j=\rho(B)$ $\left(d=\rho(C)\right)$, $\forall i,j=1,\cdots,N$, $v_r\neq\mathbf{1}$, then every solution starting in $\{x\mid x_i=1, \forall i\in \mathcal{J}, x_i\in(-1,1), \forall i\notin \mathcal{J}\}$ approaches $\frac{1}{||v_r||_{\infty}}v_r$, where $\mathcal{J}=\{j\mid v_{rj}=||v_r||_{\infty}\}$.

\item If $\sigma_i=\sigma_j=\rho(B)$ $\left(d=\rho(C)\right)$, $\forall i,j=1,\cdots,N$, $v_r\neq\mathbf{1}$, then every solution starting in $\{x\mid x_i=-1, \forall i\in \mathcal{J},\, x_i\in(-1,1), \forall i\notin \mathcal{J}\}$ approaches $-\frac{1}{||v_r||_{\infty}}v_r$.
\end{enumerate}
\end{theorem}

\begin{IEEEproof}\
First, we consider the case under Assumption~\ref{assumption1}, and the partition of the state space is shown in Fig.~\ref{figure3}.

If Assumption~\ref{assumption1} holds, then $0$ is a simple eigenvalue of $E$ with a right positive eigenvector $v_r$ when $\sigma_i=\rho(B), \forall i=1,\cdots,N$. When $\sigma_i=\sigma_j>\rho(B), \forall i,j=1,\cdots,N$, all eigenvalues of $E$ have negative real parts, and its dominant eigenvalue has the right positive eigenvector $v_r$.

\begin{figure}[htbp]
\vspace*{-3ex}
\centering
\includegraphics[width=0.48\textwidth]{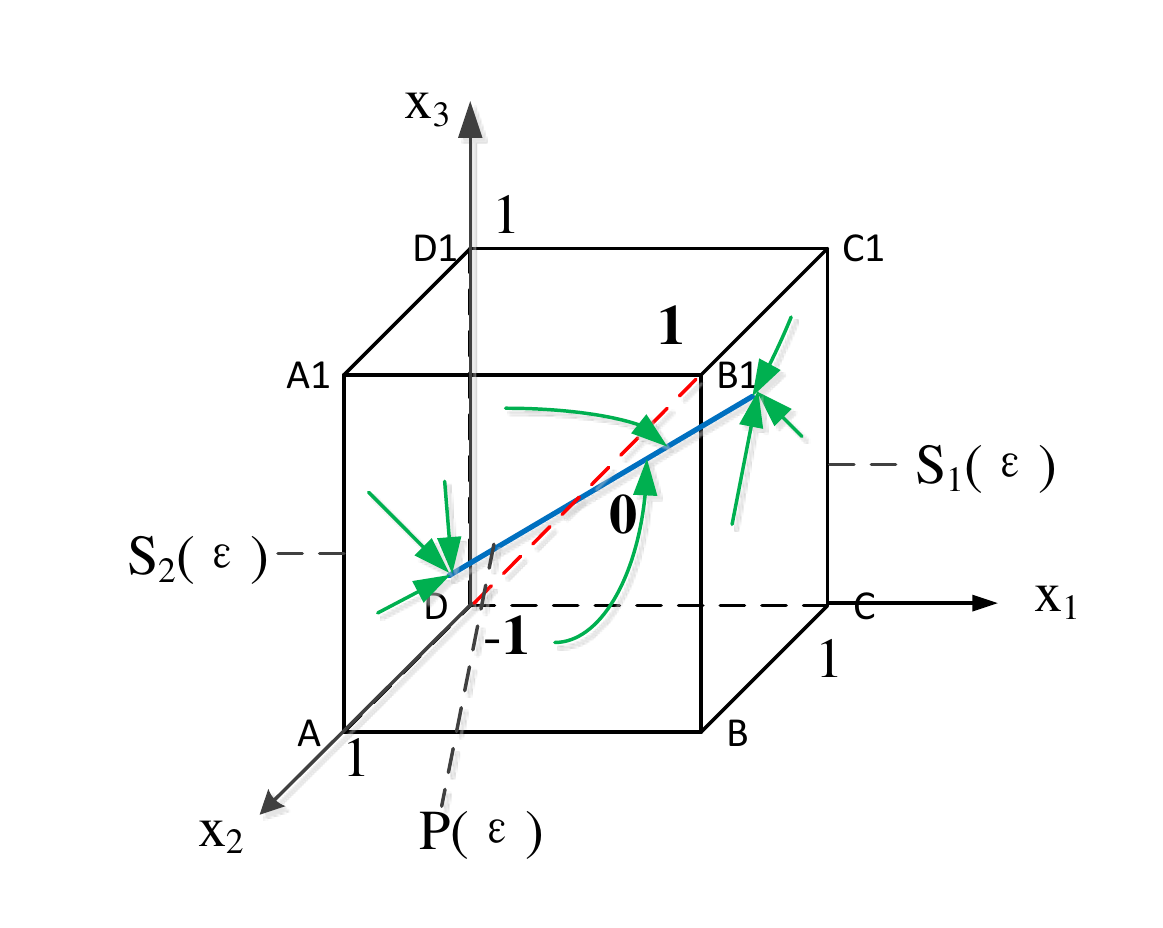}
\vspace{-3ex}
\caption{Convergence behavior of the generalized opinion dynamics (\ref{sec2:2}) with stubborn extremists in three dimensions. The blue solid line is spanned by the right positive eigenvector $v_r$, and the green curves with arrows denote several trajectories of the generalized opinion dynamics (\ref{sec2:2}).}
\label{figure3}
\end{figure}

1) Let $\mathcal{S}(\varepsilon)=\{x\mid -1+\varepsilon\leq x_i\leq1-\varepsilon\}$, where $0<\varepsilon<1$. Introduce
\begin{equation}\label{sec33:1}
  V(x)=-\frac{1}{2}\sum_{i=1}^N\gamma_i\ln(1-x_i^2).
\end{equation}
For $x\in \mathcal{S}(\varepsilon)$, $V(x)$ is positive definite, and the time derivative of $V(x)$ along (\ref{sec33:1}) is found to be
\begin{align*}
  \dot{V}&=\sum_{i=1}^N\gamma_ix_i\frac{\dot{x}_i}{1-x_i^2}\\
         &=x^{\rm T}H(\Gamma E)x.
\end{align*}
Since $H(\Gamma E)\preceq0$, we have $\dot{V}\leq 0.$ Hence, $\mathcal{S}(\varepsilon)$ is positively invariant with respect to (\ref{sec2:2}). By Lemma~\ref{lemma4}, $\dot{V}=0$ if and only if $H(\Gamma E)x=0$. Because $\rank H(\Gamma E)=\rank E$, $\dot{V}=0$ if and only if $Ex=0$. If $\sigma_i=\rho(B), \forall i=1,\cdots,N$, then $\dot{V}=0$ if and only if $x\in \mathcal{P}(\varepsilon)=\big\{x \,\big|\, x=\alpha v_r,\, \alpha\in\big[-\frac{-1+\varepsilon}{||v_r||_{\infty}},\frac{1-\varepsilon}{||v_r||_{\infty}}\big]\big\}$. If $\sigma_i=\sigma_j>\rho(B), \forall i,j=1,\cdots,N$, then $\dot{V}=0$ if and only if $x=\mathbf{0}$. By LaSalle invariance principle, the solution of (\ref{sec2:2}) with initial condition $x(0)\in \mathcal{S}(\varepsilon)$ converges into $\mathcal{P}(\varepsilon)$. If $x(0)$ is sufficiently close to $\mathbf{1}$ or $-\mathbf{1}$, then $\varepsilon\rightarrow 0$. Hence, for all $x(0)\in(-1,1)^N$, $\lim\limits_{t\rightarrow\infty}x(t)=\alpha v_r$, where $\alpha\in\big(-\frac{1}{||v_r||_{\infty}},\frac{1}{||v_r||_{\infty}}\big)$.

2) Let $\mathcal{S}_{\mathcal{I}_{a,0}}(\varepsilon)=\{x\mid |x_i|=1, \forall i\in \mathcal{I}_{a,0},\, x_i\in[\varepsilon-1,1-\varepsilon], \forall i\notin \mathcal{I}_{a,0}\}$, where $0<\varepsilon<1$. When $x(0)\in \mathcal{S}_{\mathcal{I}_{a,0}}(\varepsilon)$, the solution of (\ref{sec2:2}) satisfies $|x_i(t)|\equiv 1, \forall i\in \mathcal{I}_{a,0}$. $P$ is a permutation matrix such that
\begin{equation*}
  y=Px=\left[\begin{array}{c}y_1\\
  y_2\end{array}\right],
\end{equation*}
where $|y_{1i}|=1, \dim y_1=|\mathcal{I}_{a,0}|$. Then the network (\ref{sec2:2}) becomes
\begin{align*}
  \dot{y}
  &=\left[\begin{array}{cc}0 & 0\\
  0 & (I-\diag(y_2)^2)\end{array}\right]\left[\begin{array}{cc}\bar{E}_{11} & \bar{E}_{12}\\
  \bar{E}_{21} & \bar{E}_{22}\end{array}\right]\left[\begin{array}{c}y_1\\
  y_2\end{array}\right].
\end{align*}
Hence,
\begin{equation}\label{sec33:2}
  \dot{y}_2=(I-\diag(y_2)^2)\left(\bar{E}_{22}y_2+\bar{E}_{21}y_1^{\ast}\right),
\end{equation}
where $y_2\in[\varepsilon-1,1-\varepsilon]^{|\mathcal{I}_{a,+}|}$.

Since $\bar{E}_{22}$ is nonsingular, there exists a unique $y_2^{\ast}$ such that $\bar{E}_{22}y_2^{\ast}=-\bar{E}_{21}y_1^{\ast}$. Let $s=y_2-y_2^{\ast}$. Then we have
\begin{equation}\label{sec33:3}
  \dot{s}=(I-\diag(s+y_2^{\ast})^2)\bar{E}_{22}s,
\end{equation}
where $s\in[\varepsilon-1-y_2^{\ast},1-\varepsilon-y_2^{\ast}]^{|\mathcal{I}_{a,+}|}$. Assume that
\begin{align}\label{sec33:4}
  V(s)=-\frac{1}{2}\sum_{i=1}^{\bar{N}}\bar{\gamma}_{i+|\mathcal{I}_{a,0}|}\left[y_{2i}^{\ast}\ln\frac{\Phi_1(s_i)}{\Phi_2(s_i)}-\ln\Phi_2(s_i)\right],
\end{align}
where $\bar{N}=N-|\mathcal{I}_{a,0}|$,
$\Phi_1(s_i)=1+(s_i+y_{2i}^{\ast})$, $\Phi_2(s_i)=1-(s_i+y_{2i}^{\ast})$, and $\bar{\Gamma}=\diag\{\bar{\gamma}_1,\cdots,\bar{\gamma}_N\}=P\Gamma P^{\rm T}$. Let $\bar{\Gamma}_1=\diag\{\bar{\gamma}_1,\cdots,\bar{\gamma}_{|\mathcal{I}_{a,0}|}\}$ and $\bar{\Gamma}_2=\diag\{\bar{\gamma}_{|\mathcal{I}_{a,0}|+1},\cdots,\bar{\gamma}_N\}$.
It is easy to see from (\ref{sec33:4}) that $V(s)$ is continuously differentiable in $[\varepsilon-1-y_2^{\ast},1-\varepsilon-y_2^{\ast}]^{|\mathcal{I}_{a,+}|}$. The time derivative of $V(s)$ along (\ref{sec33:3}) is calculated as follows
\begin{align*}
  \dot{V}&=\sum_{i=1}^{N-|\mathcal{I}_{a,0}|}\bar{\gamma}_{i+|\mathcal{I}_{a,0}|}\frac{s_i\dot{s}_i}{1-(s_i+y_{2i}^{\ast})^2}\\
         &=\sum_{i=1}^{N-|\mathcal{I}_{a,0}|}s_i\bar{\gamma}_{i+|\mathcal{I}_{a,0}|}\bar{e}_{22i}s\\
         &=s^{\rm T}H(\bar{\Gamma}_2\bar{E}_{22})s,
\end{align*}
where $\bar{e}_{22i}$ is the $i$th row of matrix $\bar{E}_{22}$.
Since $H(\Gamma E)\preceq0$ and $H(\bar{\Gamma}_2\bar{E}_{22})$ is a principal submatrix of $H(\Gamma E)$, we have $H(\bar{\Gamma}_2\bar{E}_{22})\preceq0$. Hence,
\begin{align*}
  \dot{V}&=\sum_{i=1}^{N-|\mathcal{I}_{a,0}|}s_i\bar{\gamma}_{i+|\mathcal{I}_{a,0}|}\bar{e}_{22i}s\\
         &=s^{\rm T}H(\bar{\Gamma}_2\bar{E}_{22})s\\
         &\leq 0
\end{align*}
for all $s\in [\varepsilon-1-y_2^{\ast},1-\varepsilon-y_2^{\ast}]^{|\mathcal{I}_{a,+}|}$.

Thus, $[\varepsilon-1-y_2^{\ast},1-\varepsilon-y_2^{\ast}]^{|\mathcal{I}_{a,+}|}$ is positively invariant with respect to (\ref{sec33:3}). By LaSalle invariance principle, the solution of (\ref{sec33:3}) with initial condition $s(0)\in [\varepsilon-1-y_2^{\ast},1-\varepsilon-y_2^{\ast}]^{|\mathcal{I}_{a,+}|}$ converges to $y_2^{\ast}$. When $\varepsilon\rightarrow 0^{+}$, it leads to $1-\varepsilon\rightarrow1^{-}$ and $\varepsilon-1\rightarrow-1^{+}$. Hence, if $x(0)\in \{x\mid |x_i|=1, i\in \mathcal{I}_{a,0},\, x_i\in(-1,1), \forall i\notin \mathcal{I}_{a,0}\}$, then $\lim\limits_{t\rightarrow\infty}x(t)=x^{\ast}$.

3) Let $\mathcal{S}_1(\varepsilon)=\{x\mid x_j=1,\forall j\in \mathcal{J}, x_i\in [\varepsilon-1,1-\varepsilon], \forall i\notin \mathcal{J}\}$, where $0<\varepsilon<1$. For $x\in \mathcal{S}_1(\varepsilon)$, there exists a permutation matrix $P$ such that the network (\ref{sec2:2}) has the following form
\begin{align*}
  \dot{y}
  &=\left[\begin{array}{cc}0 & 0\\
  0 & (I-\diag(y_2)^2)\end{array}\right]\left[\begin{array}{cc}\bar{E}_{11} & \bar{E}_{12}\\
  \bar{E}_{21} & \bar{E}_{22}\end{array}\right]\left[\begin{array}{c}y_1\\
  y_2\end{array}\right],
\end{align*}
where
\begin{equation*}
  y=Px=\left[\begin{array}{c}y_1\\
  y_2\end{array}\right]
\end{equation*}
with $y_1=\mathbf{1}, \dim y_1=|\mathcal{J}|$. Then we can get
\begin{equation}\label{sec33:5}
  \dot{y}_2=(I-\diag(y_2)^2)\left(\bar{E}_{22}y_2+\bar{E}_{21}\mathbf{1}\right).
\end{equation}
Because $v_r\neq\mathbf{1}$ and $\mathcal{J}=\{j\mid v_{rj}=||v_r||_{\infty}\}$, there exists a unique solution $y_2^{\ast}$ such that $\bar{E}_{22}y_2^{\ast}+\bar{E}_{21}\mathbf{1}=0$. From a geometric point of view, the point
\begin{equation*}
  x^{\ast}=P^{\rm T}\left[\begin{array}{c}\mathbf{1}\\
  y_2^{\ast}\end{array}\right]
\end{equation*}
is the intersection point of the line $\alpha v_r$ with $(-1,1]^N$. Similar to the arguments presented in 2), we can obtain that every solution starting in $\mathcal{S}_1(\varepsilon)$ approaches $x^{\ast}=\frac{1}{||v_r||_{\infty}}v_r$. When $\varepsilon\rightarrow 0^{+}$, we can have that every solution starting in $\{x\mid x_i=1, \forall i\in \mathcal{J}, x_i\in(-1,1), \forall i\notin \mathcal{J}\}$ approaches $\frac{1}{||v_r||_{\infty}}v_r$.

4) Let $\mathcal{S}_2(\varepsilon)=\{x\mid x_i=-1,\forall i\in \mathcal{J}, x_i\in [\varepsilon-1,1-\varepsilon], \forall i\notin \mathcal{J}\}$, where $0<\varepsilon<1$. For $x(0)\in \mathcal{S}_2(\varepsilon)$, $x_i\equiv-1, \forall i\in \mathcal{J}$. There exists a permutation matrix $P$ such that the network (\ref{sec2:2}) takes the following form
\begin{align*}
  \dot{y}
  &=\left[\begin{array}{cc}0 & 0\\
  0 & (I-\diag(y_2)^2)\end{array}\right]\left[\begin{array}{cc}\bar{E}_{11} & \bar{E}_{12}\\
  \bar{E}_{21} & \bar{E}_{22}\end{array}\right]\left[\begin{array}{c}y_1\\
  y_2\end{array}\right],
\end{align*}
where
\begin{equation*}
  y=Px=\left[\begin{array}{c}y_1\\
  y_2\end{array}\right]
\end{equation*}
with $y_1=-\mathbf{1}, \dim y_1=|\mathcal{J}|$. Hence,
\begin{equation}\label{sec33:6}
  \dot{y}_2=(I-\diag(y_2)^2)\left(\bar{E}_{22}y_2-\bar{E}_{21}\mathbf{1}\right).
\end{equation}
Because $\sigma_i=\sigma_j=\rho(B)$, $\forall i,j=1,\cdots,N$, $v_r\neq\mathbf{1}$, there exists a unique solution $y_2^{\ast}$ such that $\bar{E}_{22}y_2^{\ast}-\bar{E}_{21}\mathbf{1}=0$. From a geometric viewpoint, the point
\begin{equation*}
  x^{\ast}=P^{\rm T}\left[\begin{array}{c}-\mathbf{1}\\
  y_2^{\ast}\end{array}\right]
\end{equation*}
is the intersection point of the line $\alpha v_r$ with $[-1,1)^N$. Similar to 3), we can derive that every solution starting in $\{x\mid x_i=-1, \forall i\in \mathcal{J}, x_i\in(-1,1), \forall i\notin \mathcal{J}\}$ approaches $-\frac{1}{||v_r||_{\infty}}v_r$.

Second, we deal with the case that Assumption~\ref{assumption2} holds. Similar to Theorem~\ref{theorem1}, we can know that $E=C-dI$ has the eigenvalue $\rho(C)-d$ with eigenvector $v_r$. If $d=\rho(C)$, then all eigenvalues of $E$ have non-positive real parts, and $0$ is a simple eigenvalue of $E$ with a right positive eigenvector $v_r$. If $d>\rho(C)$, then all eigenvalues of $E$ have negative real parts, and its dominant eigenvalue has a right positive eigenvector $v_r$.

Finally, we can adopt the similar methods presented above to get 1), 2), 3) and 4) under Assumption~\ref{assumption2}, and we omit details here. This completes the proof of Theorem~\ref{theorem3}.
\end{IEEEproof}

When the signed graph $\mathcal{G}(B)$ is undirected, the next corollary is an immediate result from Theorem~\ref{theorem3}. The proof is skipped here since it is similar to that of Corollary \ref{corollary1}.

\begin{corollary}\label{corollary3}
Suppose that Assumption~\ref{assumption1} (Assumption~\ref{assumption2}) holds and the signed graph $\mathcal{G}(B)$ is undirected.
\begin{enumerate}[1)]
\item If $x(0)\in (-1,1)^N$, then $\lim\limits_{t\rightarrow\infty}x(t)=\alpha v_r$, where $v_r$ is the right positive eigenvector of $B$ $(C)$, and $\alpha\in\big(-\frac{1}{||v_r||_{\infty}},\frac{1}{||v_r||_{\infty}}\big)$.

\item If $\mathcal{I}_{a,0}\neq\varnothing$, $\mathcal{I}_{a,+}\neq\varnothing$, there exists a permutation matrix $P$ such that the network (\ref{sec2:2}) has the form of (\ref{sec3:2}), $\bar{E}_{22}$ is nonsingular, and $-\mathbf{1}<-\bar{E}_{22}^{-1}\bar{E}_{21}\mathbf{1}<\mathbf{1}$, then $\lim\limits_{t\rightarrow\infty}x(t)=x^{\ast}$ for all $x(0)\in \{x\mid |x_i|=1, i\in \mathcal{I}_{a,0}, x_i\in(-1,1), \forall i\notin \mathcal{I}_{a,0}\}$, where $x^{\ast}=P^{\rm T}y^{\ast}$, and $y^{\ast}$ is the equilibrium point of (\ref{sec3:2}).

\item If $\sigma_i=\sigma_j=\rho(B)$ $\left(d=\rho(C)\right)$, $\forall i,j=1,\cdots,N$, $v_r\neq\mathbf{1}$, then every solution starting in $\{x\mid x_i=1, \forall i\in \mathcal{J}, x_i\in(-1,1), \forall i\notin \mathcal{J}\}$ approaches $\frac{1}{||v_r||_{\infty}}v_r$, where $\mathcal{J}=\{j\mid v_{rj}=||v_r||_{\infty}\}$.

\item If $\sigma_i=\sigma_j=\rho(B)$ $\left(d=\rho(C)\right)$, $\forall i,j=1,\cdots,N$, $v_r\neq\mathbf{1}$, then every solution starting in $\{x\mid x_i=-1, \forall i\in \mathcal{J}, x_i\in(-1,1), \forall i\notin \mathcal{J}\}$ approaches $-\frac{1}{||v_r||_{\infty}}v_r$.
\end{enumerate}
\end{corollary}

The stubborn extremists scenario illustrates the circumstance that the extreme opinions cannot be persuaded by other opinions, while the other opinions are readily changed. This scenario is adequate when there exists two competitive extreme opinions, such as two opposite political parties. In this scenario it is demonstrated in Theorem~\ref{theorem3} and Corollary~\ref{corollary3} that the generalized opinion dynamics (\ref{sec2:2}) has the assured convergence under different initial conditions.

\subsection{Further Discussion}\label{section34}

In this subsection, we will discuss several important aspects of the results derived in the preceding subsections.

\begin{remark}\label{remark1}
For different state-dependent susceptibility, the dynamical behaviors of the generalized opinion model (\ref{sec2:2}) are very different. It can be seen from the proofs of Theorem~\ref{theorem1}--Theorem~\ref{theorem3} that the number of equilibrium points is different for different state-dependent susceptibility. Specifically, the number of equilibrium points is smallest in the stubborn neutrals scenario, while the number of equilibrium points is biggest in the stubborn extremists scenario. As a result, distinctive Lyapunov functions (compare (\ref{sec31:1}), (\ref{sec31:4}), (\ref{sec32:1}), (\ref{sec32:2}), (\ref{sec32:3}), (\ref{sec33:1}), (\ref{sec33:4})) are constructed to effectively handle these different scenarios with a view to proving the convergence of the generalized opinion dynamics (\ref{sec2:2}).
\end{remark}

\begin{remark}\label{remark2}
If for each $x_0\in\mathbb{R}^{n}$, there exists $T>0$ such that $x\in\mathbb{R}^{n}_{+}$ or $x\in\mathbb{R}^{n}_{-}$ $\forall t\geq T$, then the system is said to achieve an unanimous opinion. Our Theorem~\ref{theorem1}--Theorem~\ref{theorem3} reveal that the network (\ref{sec2:2}) achieves a unanimous opinion. Moreover, although all agents achieve a unanimous opinion, the degree of opinion for each agent is different when the right positive eigenvector $v_r\neq \mathbf{1}$.
\end{remark}

\begin{remark}\label{remark3}
In Theorem~\ref{theorem1}--Theorem~\ref{theorem3}, one may use the condition ``$E$ is real orthogonally diagonalizable" to replace the condition ``$H(\Gamma E)\preceq0$ and $\rank H(\Gamma E)=\rank E$." Because $E$ is real orthogonally diagonalizable, there exists a real orthogonal matrix $U$ such that $UEU^{\rm T}$ is a diagonal matrix. By Lemma~\ref{lemma5} and Lemma~\ref{lemma6}, we have that $H(E)\preceq0$ and $\rank H(E)=\rank E$.
\end{remark}

\begin{remark}\label{remark4}
If the signed graph $\mathcal{G}(B)$ is undirected, then matrix $E$ is a real symmetric matrix and orthogonally diagonalizable. If the signed graph $\mathcal{G}(B)$ is directed and $B$ is real orthogonally diagonalizable, then  matrix $E$ is orthogonally diagonalizable.
\end{remark}

\begin{remark}\label{remark5}
If matrix $E$ is not symmetric, then it is not easy to verify whether the matrix is real orthogonally diagonalizable. However, by using the Matlab LMI (linear matrix inequality) toolbox, it is straightforward to find a positive diagonal matrix $\Gamma$ such that $H(\Gamma E)\preceq0$ and $\rank H(E)=\rank E$. See Example~\ref{exm2} in next section for details.
\end{remark}

\begin{remark}\label{remark6}
For a special class of matrices, there exists a positive diagonal matrix $\Gamma$ such that $H(\Gamma E)\preceq0$. If $\sigma_i=\sigma_j>\rho(B)$, $\forall i,j=1,\cdots,N$, and $B\geq 0$, then $E=B-\Sigma$ is a nonsingular M-matrix. It is known from the property of M-matrices \cite{olesky2009m} that there exists a positive diagonal matrix $\Gamma$ such that $H(\Gamma E)\prec0$. Moreover, $\rank H(E)=\rank E$.
\end{remark}

\begin{remark}\label{remark7}
Unlike the model in \cite{altafini2015predictable}, our model requires that all eigenvalues of matrix $E$ have non-positive real parts. If there exists an eigenvalue of $E$ which has the positive real part, then the evolution of one agent may approach infinity as time goes to infinity. This will contradict the assumption in the stubborn neutrals scenario $A(x)=\diag(x)^2$. Moreover, in practice, the opinion about a specific matter should be measured as a finite number rather than infinity.
\end{remark}

\begin{remark}\label{remark8}
The major difference between our work and \cite{amelkin2017polar} is that the work in this paper can deal with antagonistic interactions in social networks. Moreover, the model in this paper is more general than the one in \cite{amelkin2017polar}. For example, all eigenvalues of matrix $E$ may have negative real parts in our model. However, when all interactions among agents are cooperative, our results may be more conservative than those in \cite{amelkin2017polar}, especially in the case of directed graphs.
\end{remark}

\section{Numerical Examples}\label{section4}

This section will provide two examples to illustrate the effectiveness of the theoretical results. The first example comes from the real-world, that is, Zachary's Karate Club \cite{zachary1977information}, while the second example is taken from \cite{altafini2015predictable}.

\begin{example}\label{exm1}
In \cite{zachary1977information}, Zachary studied the relationships in a karate club, and obtained some results about fission. The club consists of 34 members. There are two important persons in the club: one is the club president John A., and the other is Mr. Hi who is a part-time karate instructor. There is an conflict between John A., and Mr.~Hi over the price of karate lessons. As the instructor, Mr.~Hi wishes to raise prices. However, as the club's chief administrator, John A. wishes to stabilize prices. Assume that the viewpoint of John A. is $1$, the viewpoint of Mr.~Hi is $-1$, and the viewpoint of other members is included in $[-1,1]$. In \cite{amelkin2017polar}, the three scenarios were studied, that is, stubborn positives scenario $A(x)=0.5(I-\diag(x))$, stubborn neutrals scenario $A(x)=\diag(x)^2$, and stubborn extremists scenario $A(x)=(I-\diag(x)^2)$. However, antagonistic interactions among 34 members are also likely, but were not considered in \cite{amelkin2017polar}. In this example, we assume that there exist three pairs of antagonistic interactions, that is, Members 1 and 2, Members 1 and 32, and Members 33 and 34. We examine three different scenarios, that is, stubborn positives scenario, stubborn neutrals scenario, and stubborn extremists scenario.

Fig.~\ref{figure4} shows the evolution of each person's viewpoint in the stubborn positives scenario, and illustrates the results of Theorem~\ref{theorem1}. It is clear from Fig.~\ref{figure4} that the evolution of each person's viewpoint achieves the same viewpoint with different levels of the opinion. Fig.~\ref{figure5} depicts the evolution of each person's viewpoint in the stubborn neutrals scenario, exemplifying the results of Theorem~\ref{theorem2}. When all persons' initial viewpoints are the same as John A's, that is, $x_i(0)>0$, the evolution of each person's viewpoint reaches the same viewpoint as John A but at different levels of the opinion. Furthermore, as an illustration of the results of Theorem~\ref{theorem3}, Fig.~\ref{figure6} displays the evolution of each person's viewpoint in the stubborn extremists scenario. In particular, when there exist $i,j$ such that $x_i(0)=1, x_j(0)=-1$, the evolutions of these persons' viewpoint are unable to attain the same opinion, as seen from Fig.~\ref{figure6}(d). Note that this phenomenon cannot be deduced from Theorem~\ref{theorem3}, which is left for future research.

\begin{figure}[htbp]
\centering
\subfigure[]{\includegraphics[width=0.42\textwidth]{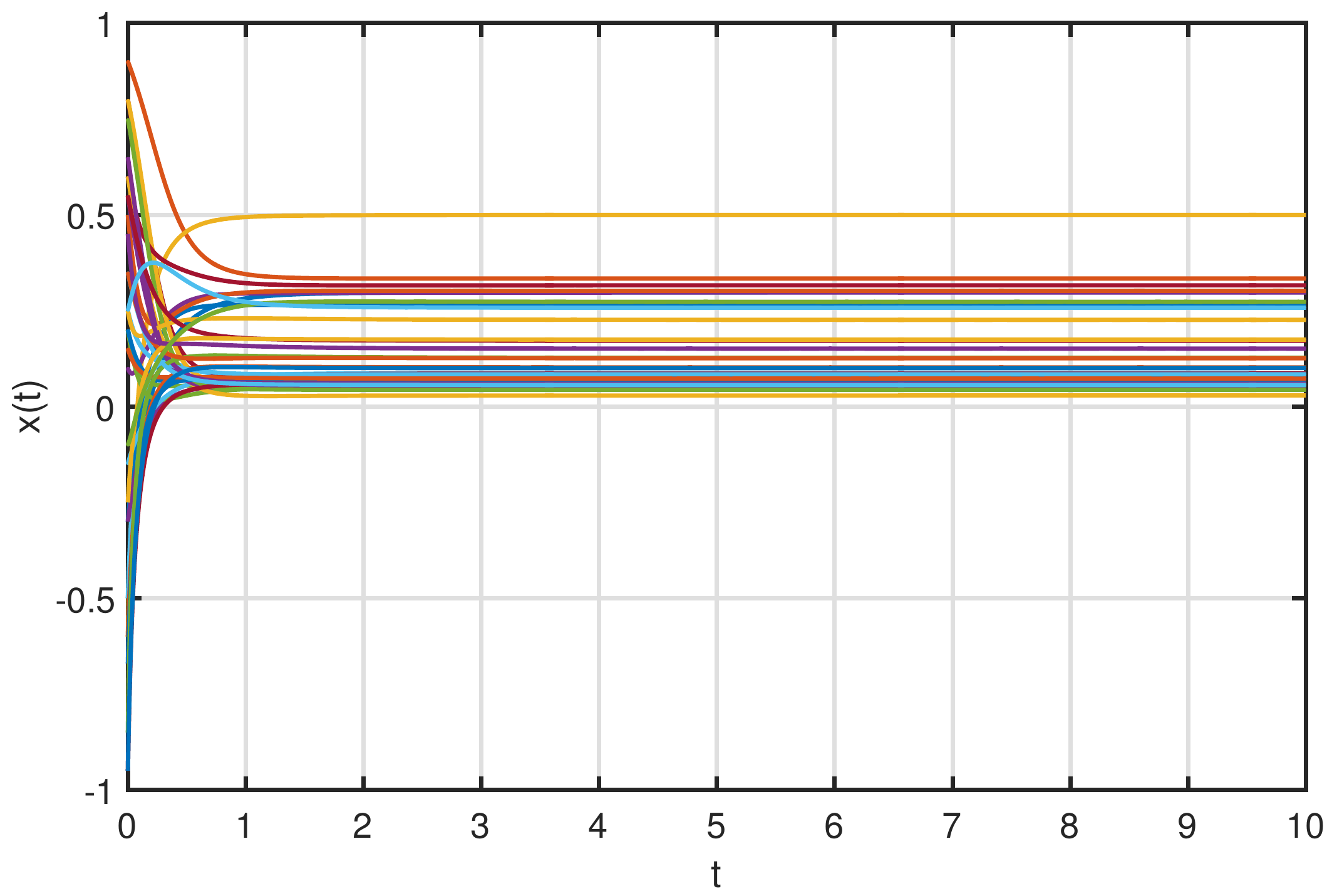}}
\subfigure[]{\includegraphics[width=0.42\textwidth]{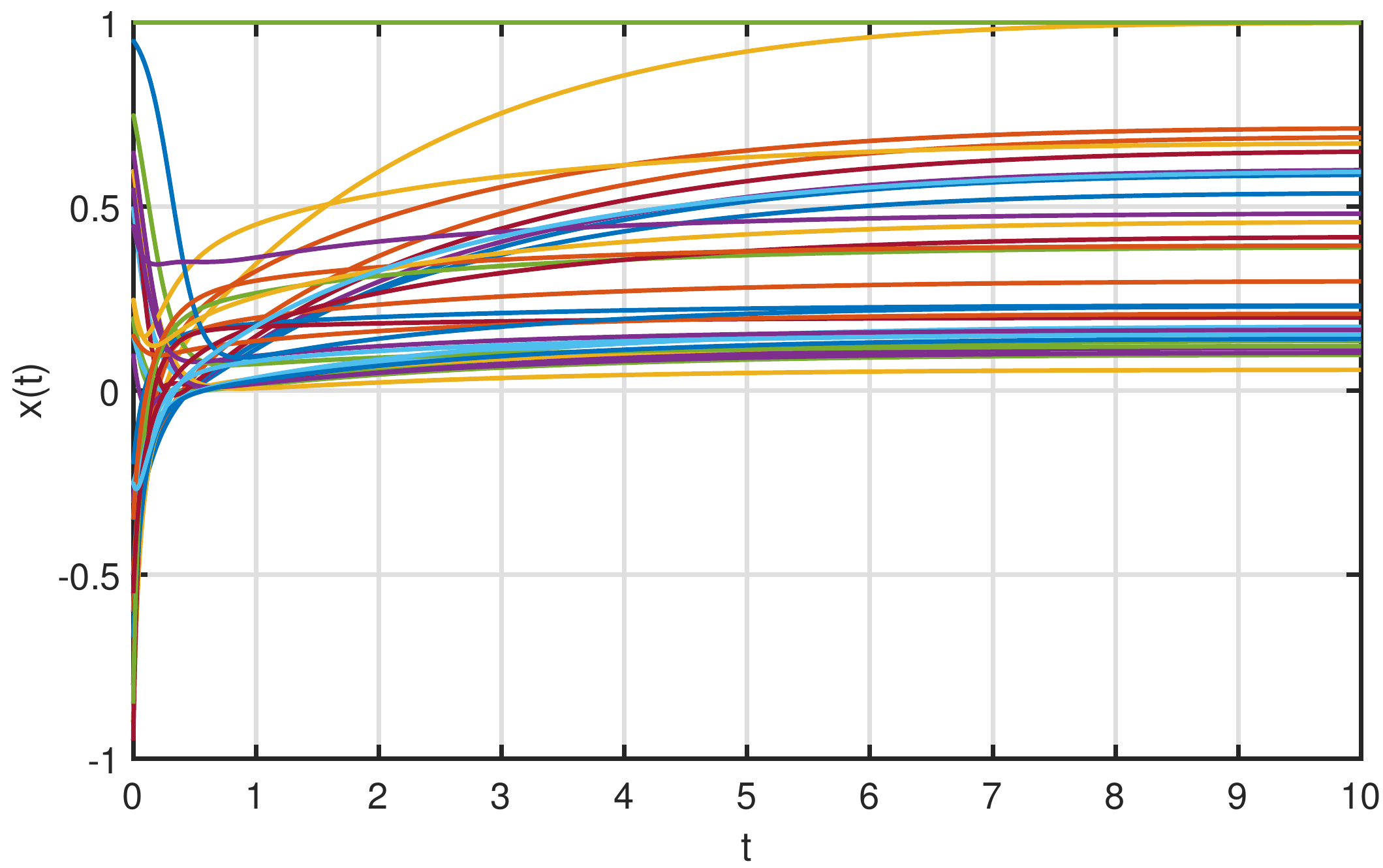}}
\vspace{-1ex}
\caption{Stubborn positives scenario in Example~\ref{exm1}: (a) Evolution of each person's viewpoint for all $|x_i(0)|<1$; (b) Evolution of each person's viewpoint under the condition that $x_{33}(0)=1, |x_i(0)|<1, \forall i\neq33$.}
\label{figure4}
\vspace{-1ex}
\end{figure}

\begin{figure}[htbp]
\centering
\subfigure[]{\includegraphics[width=0.42\textwidth]{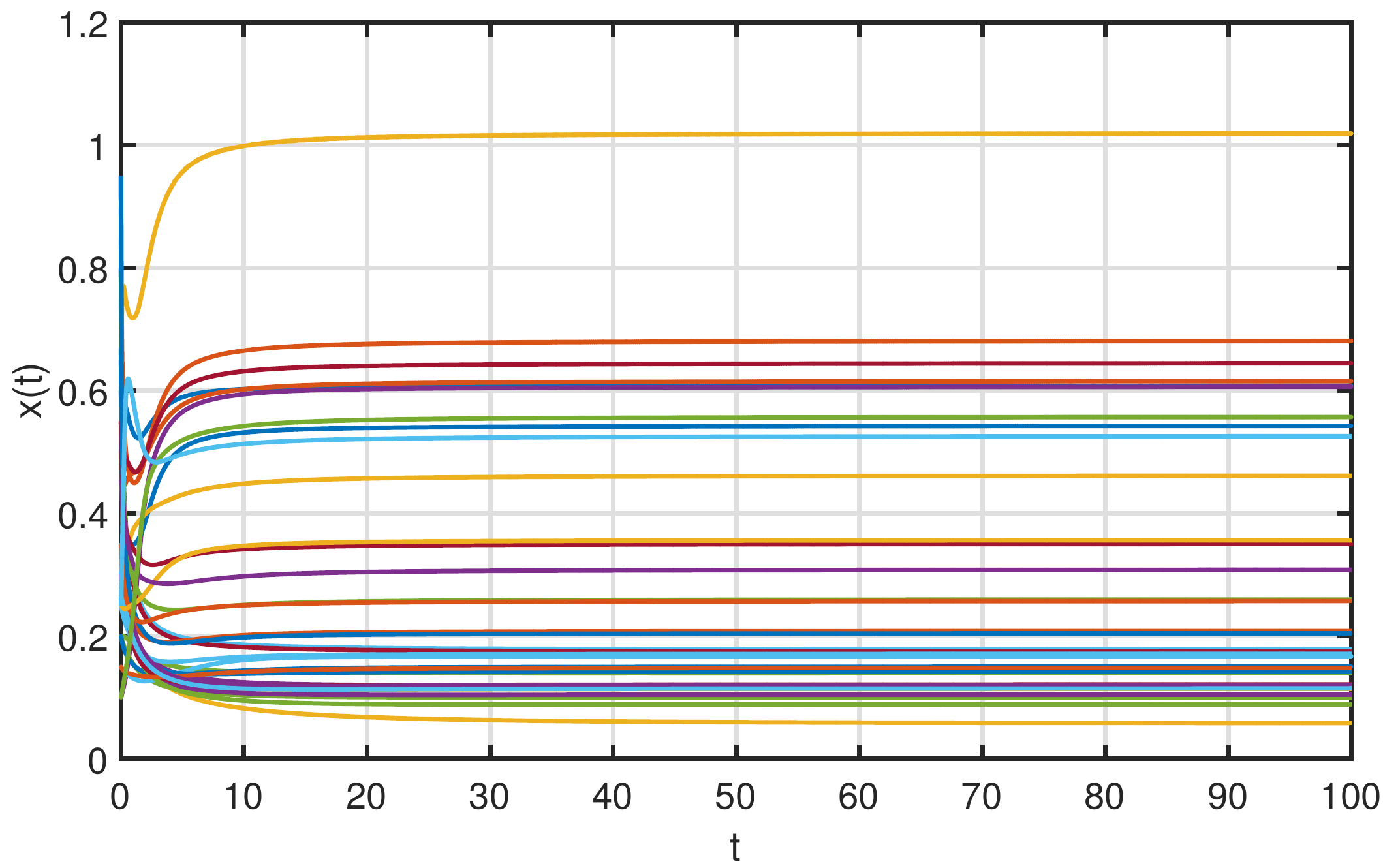}}
\subfigure[]{\includegraphics[width=0.42\textwidth]{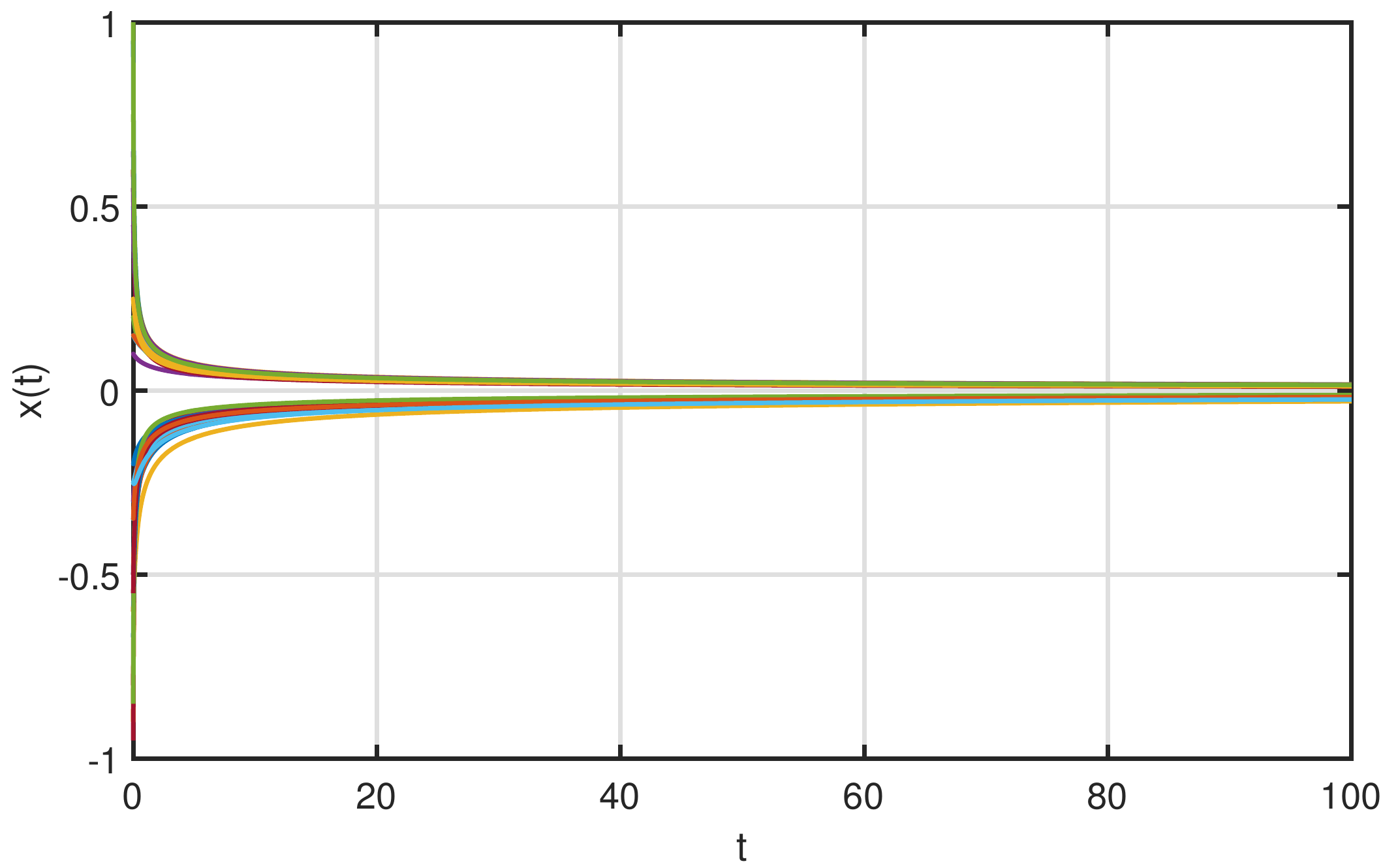}}
\vspace{-1ex}
\caption{Stubborn neutrals scenario in Example~\ref{exm1}: (a) Evolution of each person's viewpoint for all $x_i(0)>0$; (b) Evolution of each person's viewpoint under the condition that $\exists i,j~\text{s.t.}~x_i(0)>0, x_j(0)<0$.}
\label{figure5}
\end{figure}

\begin{figure}[htbp]
\centering
\subfigure[]{\includegraphics[width=0.42\textwidth]{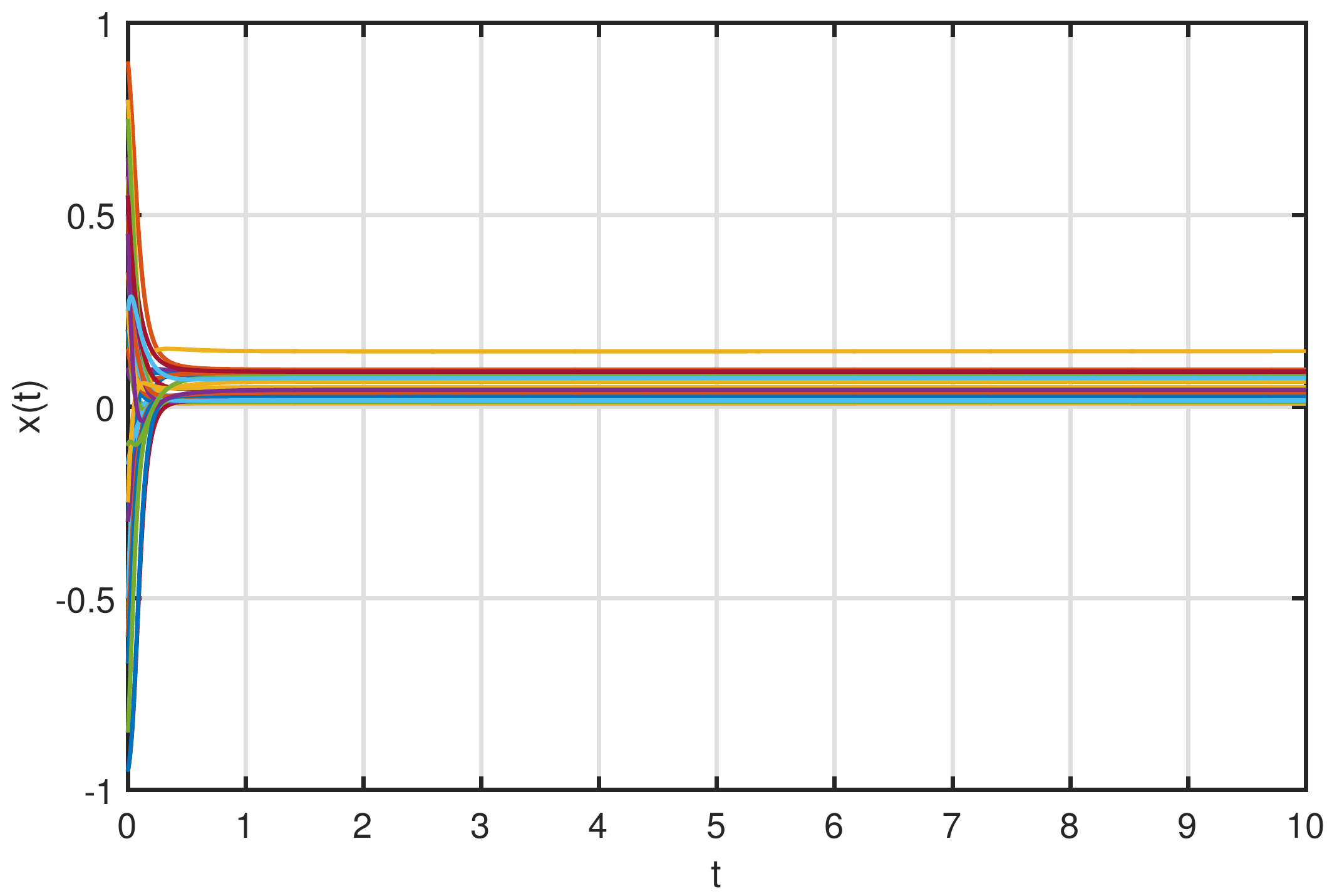}}
\subfigure[]{\includegraphics[width=0.42\textwidth]{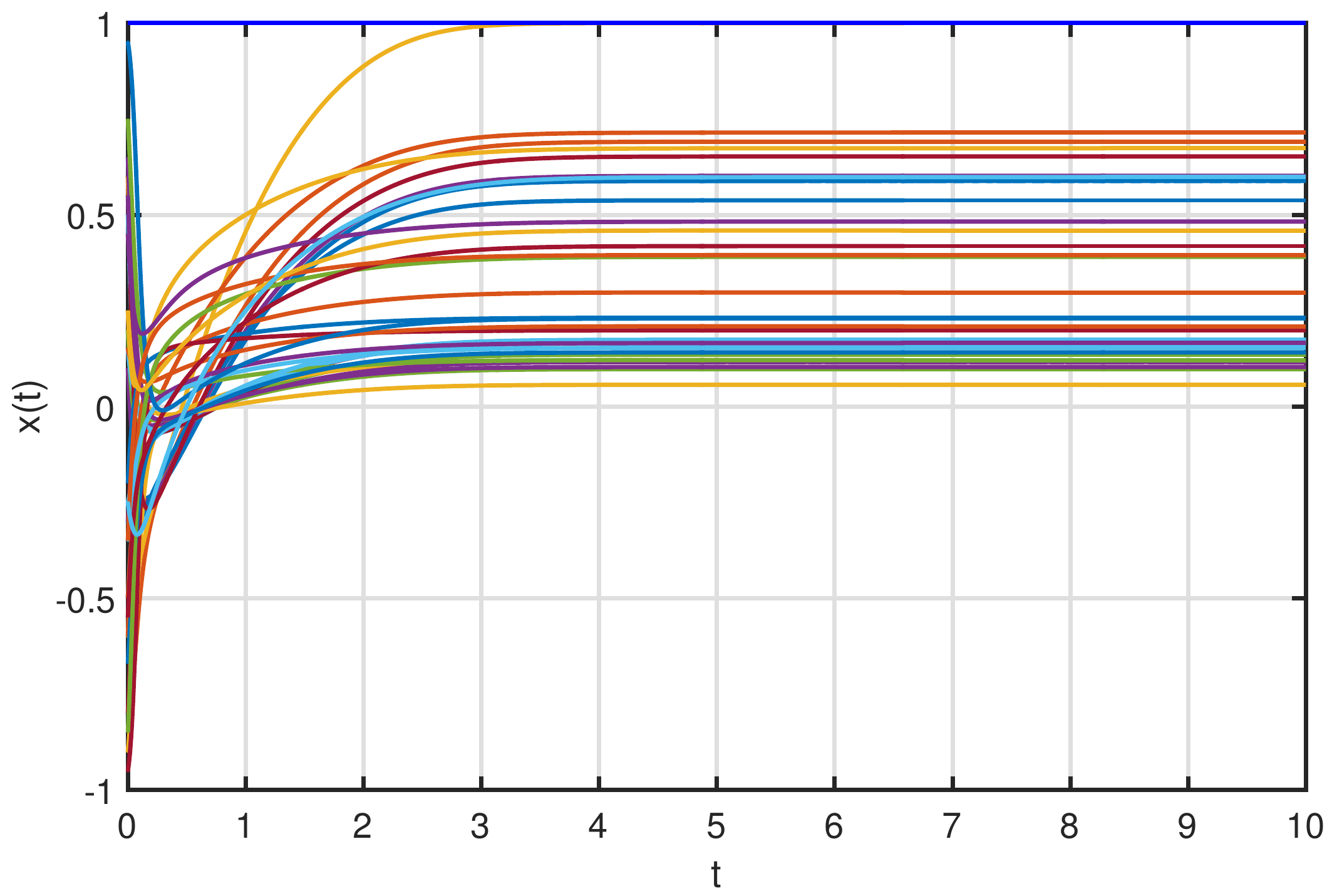}}
\subfigure[]{\includegraphics[width=0.42\textwidth]{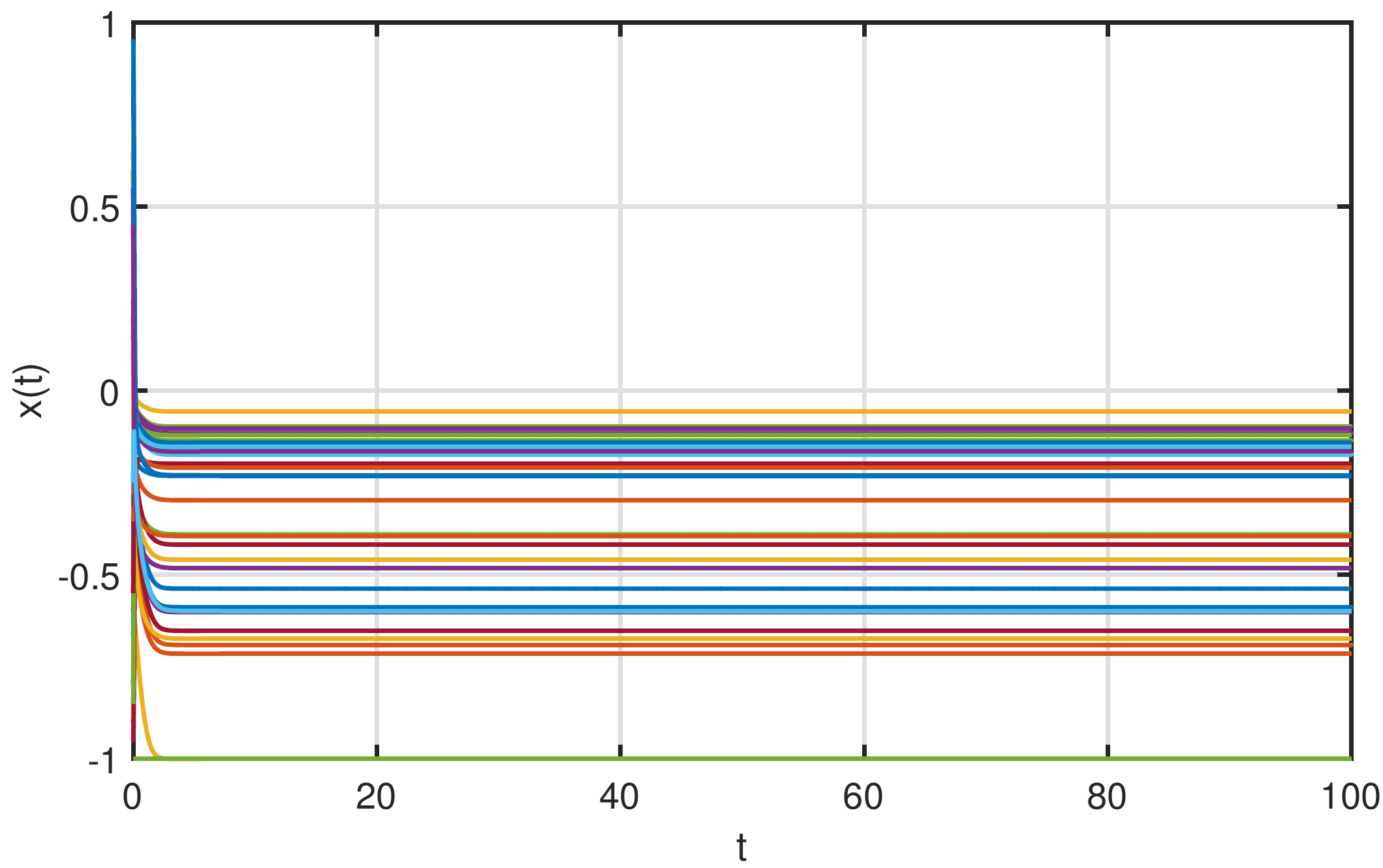}}
\subfigure[]{\includegraphics[width=0.42\textwidth]{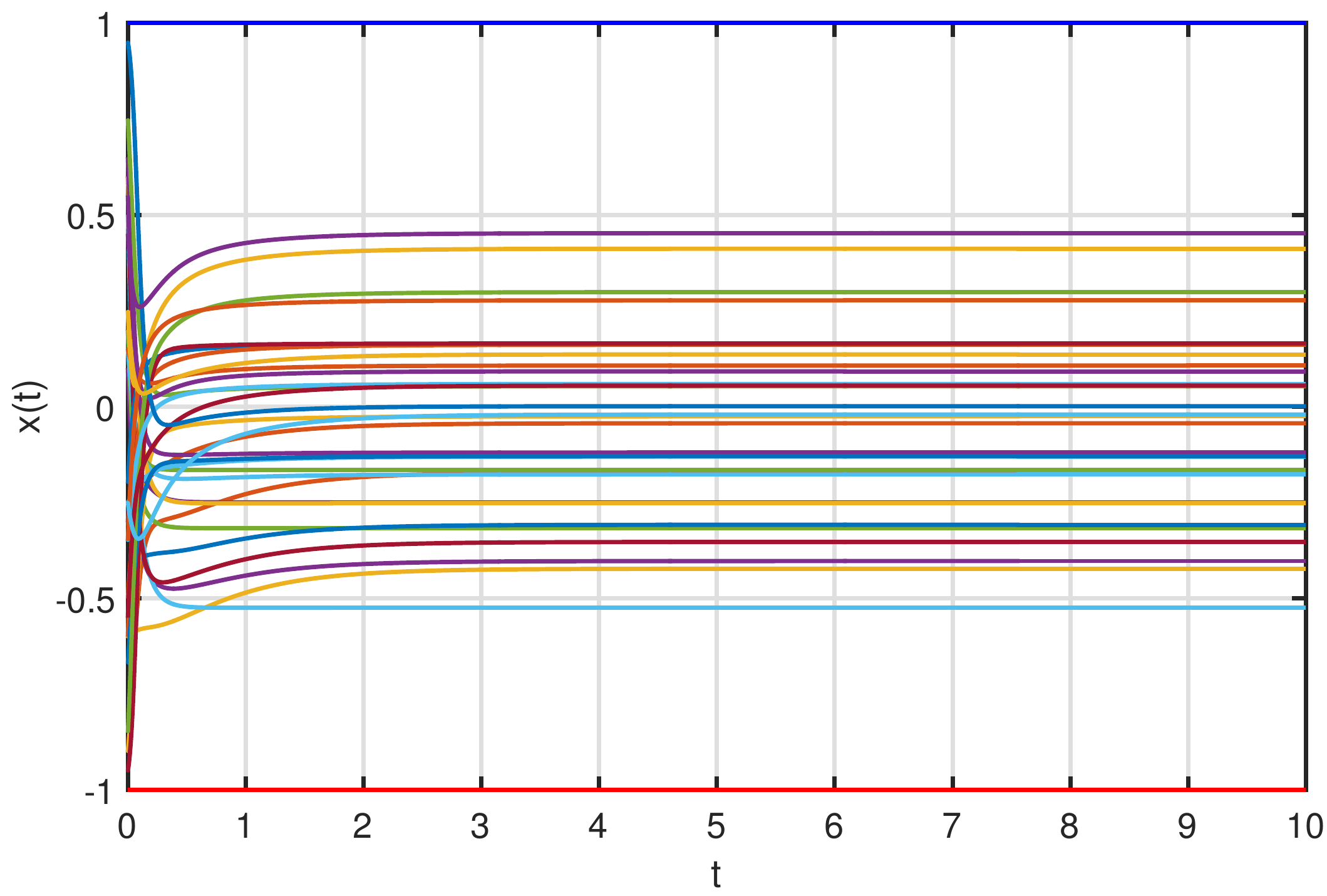}}
\vspace{-1ex}
\caption{Stubborn extremists scenario in Example~\ref{exm1}: (a) Evolution of each person's viewpoint for all $|x_i(0)|<1$; (b) Evolution of each person's viewpoint under the condition that $x_{33}(0)=1, |x_i(0)|<1, \forall i\neq33$; (c) Evolution of each person's viewpoint under the condition that $x_{33}(0)=-1, |x_i(0)|<1, \forall i\neq33$; (d) Evolution of each person's viewpoint under the condition that $\exists i,j~\text{s.t.}~x_i(0)=1, x_j(0)=-1$.}
\label{figure6}
\end{figure}
\end{example}

\begin{example}\label{exm2}
Consider the opinion dynamics
\begin{equation}\label{sec4:1}
  \dot{x}=A(x)(-\Sigma+B)x,
\end{equation}
where $B$ is the adjacency matrix of signed digraph $\mathcal{G}(B)$, and
\begin{equation*}
  B=\left[\begin{array}{ccc}
          0 & 1.7877 & -0.6743\\
          -0.7678 & 0 & 0.7354\\
          0.5878 & 0 & 0\end{array}\right].
\end{equation*}
As pointed out in \cite{altafini2015predictable}, $B$ is not eventually nonnegative. But one can choose $D=\diag\{0.2688,1.002,1.3272\}$ such that $C=B+D$ is eventually positive. The matrix $C$ has the spectral radius $\rho(C)=1.5817$ and the corresponding right positive eigenvector $v_c=[0.3350\;0.5378\;0.7737]$ with $||v_c||_{\infty}=1$. If we select  $\Sigma=1.5817I-D$, then
\begin{equation*}
  E=\left[\begin{array}{ccc}
  -1.3129 & 1.7877 & -0.6743\\
  -0.7678 & -0.5797 &  0.7354\\
  0.5878  &  0  & -0.2545\end{array}\right]\preceq0,
\end{equation*}
and $0$ is a simple positive eigenvalue of $E=-\Sigma+B$ with the corresponding right positive eigenvector being $v_c$. We can find a positive diagonal matrix $\Gamma=\diag\{4.2681,8.1972,11.5733\}$ such that $H(\Gamma E)\preceq0$ and $\rank H(\Gamma E)=\rank E$. Fig.~\ref{figure7}--Fig.~\ref{figure9} plot the evolution curves of each agent under the stubborn positives scenario $A(x)=0.5(I-\diag(x))$, the stubborn neutrals scenario $A(x)=\diag(x)^2$, and the stubborn extremists scenario $A(x)=(I-\diag(x)^2)$, respectively. So Fig.~\ref{figure7}--Fig.~\ref{figure9} illustrate the results of Theorem~\ref{theorem1}--Theorem~\ref{theorem3}, respectively.

It is interesting to note that system (\ref{sec4:1}) exhibits some new behaviors when compared with the following linear model given in \cite{altafini2015predictable}:
\begin{equation}\label{sec4:2}
  \dot{x}=(-\Sigma+B)x.
\end{equation}
For example, as pointed out in \cite[Example~2]{altafini2015predictable}, the states of (\ref{sec4:2}) hold to the orthant pair $\mathbb{R}^{3}_{-,+}$. However, when $A(x)=\diag(x)^2$, the states of (\ref{sec4:1}) hold to the orthant $\mathbb{R}^{3}_{+}$ for initial condition $x_i(0)>0$ and to the orthant $\mathbb{R}^{3}_{-}$ for initial condition $x_i(0)<0$ (see Fig.~\ref{figure8}). Moreover, when $A(x)=(I-\diag(x)^2)$, every solution starting in $\{x\mid x_i=1, \forall i\in \mathcal{J}, x_i\in(-1,1), \forall i\notin \mathcal{J}\}$ approaches $v_c$, where $\mathcal{J}=\{3\}$ (see Fig.~\ref{figure9}(b)). And every solution starting in $\{x\mid x_i=-1, \forall i\in \mathcal{J}, x_i\in(-1,1), \forall i\notin \mathcal{J}\}$ approaches $-v_c$ (see Fig.~\ref{figure9}(c)).

\begin{figure}[htbp]
\centering
\subfigure[]{\includegraphics[width=0.42\textwidth]{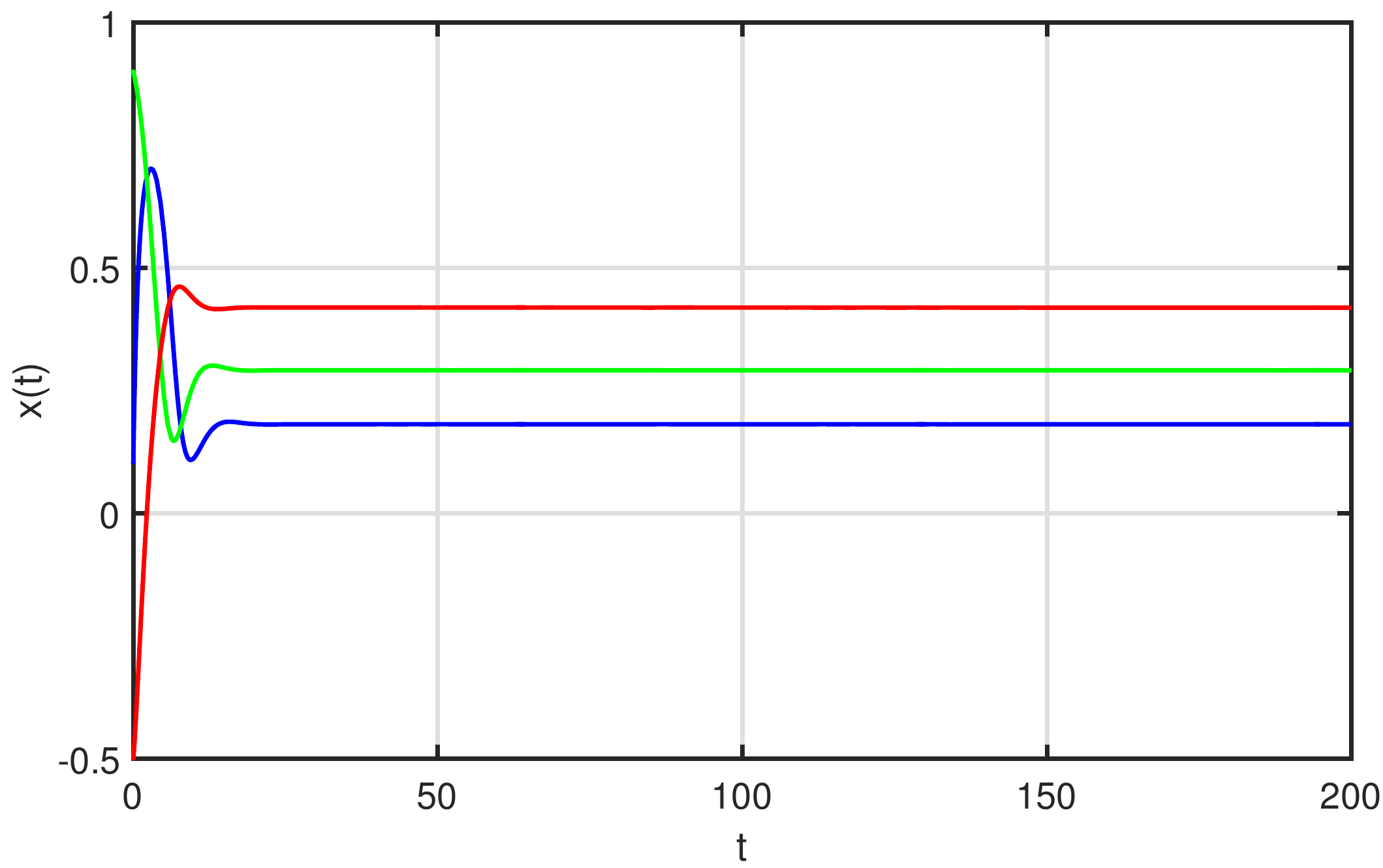}}
\subfigure[]{\includegraphics[width=0.42\textwidth]{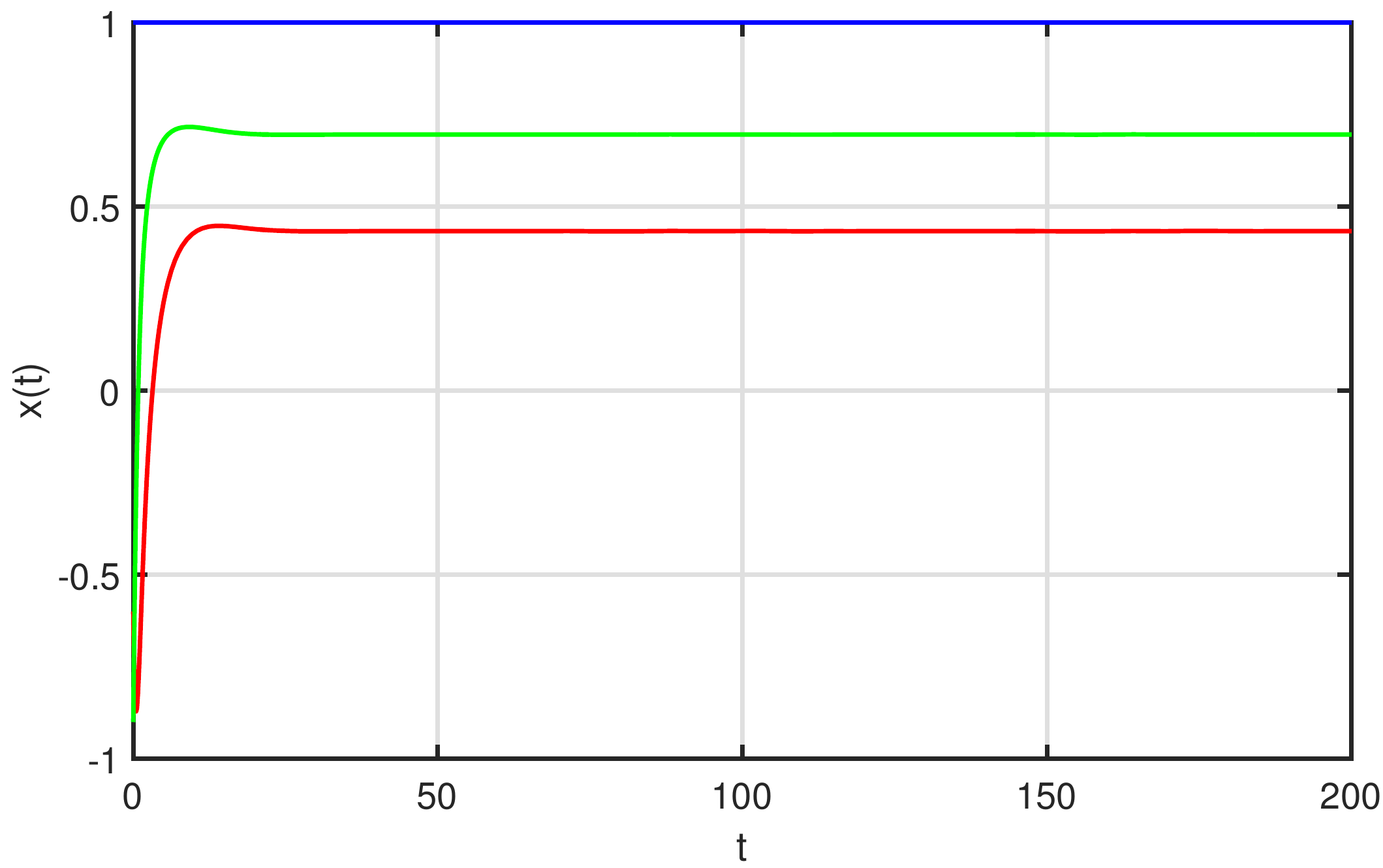}}
\caption{Stubborn positives scenario in Example~\ref{exm2}: (a) Evolution of agent for all $|x_i(0)|<1$; (b) Evolution of each agent under the condition that $x_3(0)=1, |x_i(0)|<1, \forall i\neq3$.}
\label{figure7}
\end{figure}

\begin{figure}[htbp]
\centering
\subfigure[]{\includegraphics[width=0.42\textwidth]{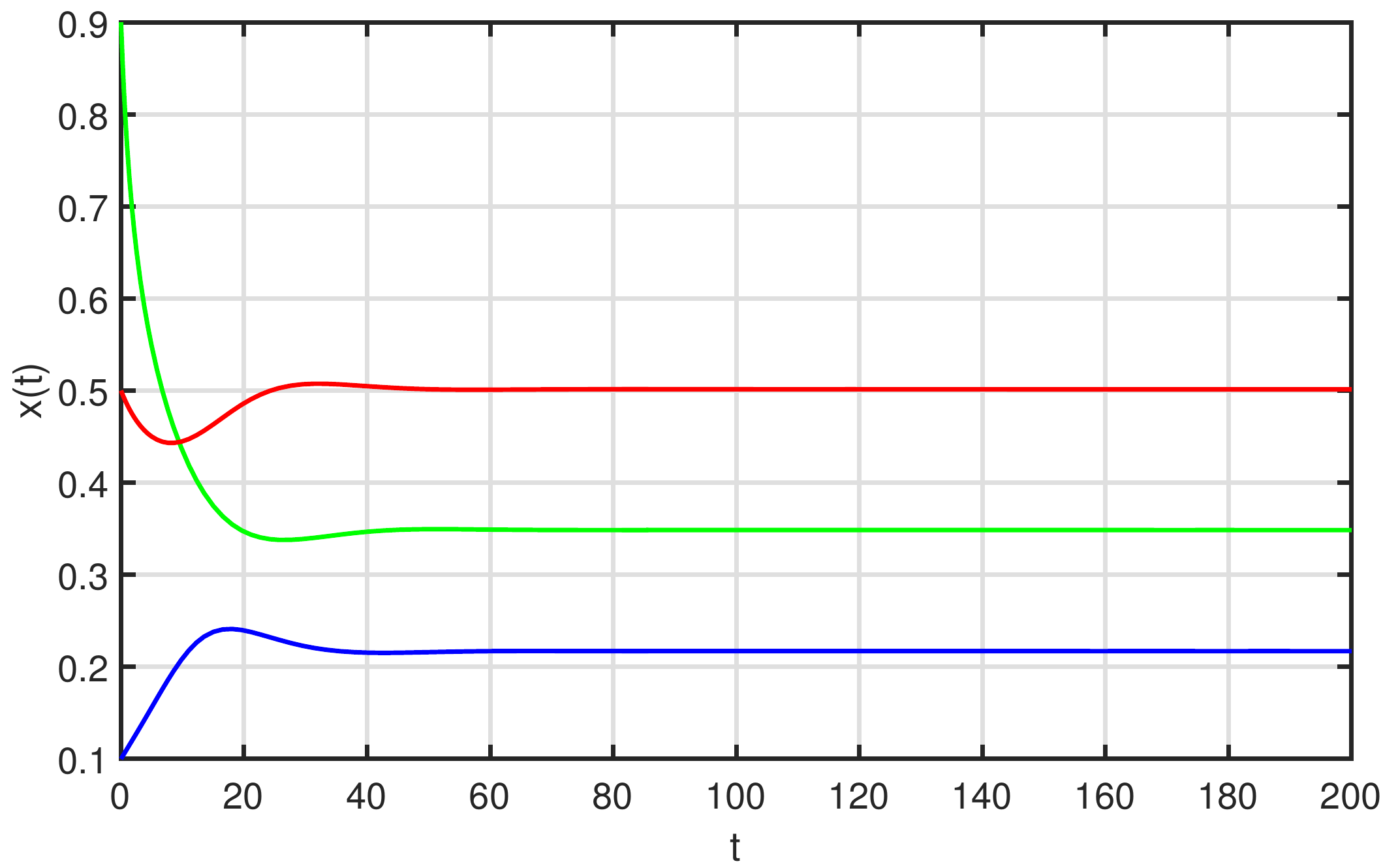}}
\subfigure[]{\includegraphics[width=0.42\textwidth]{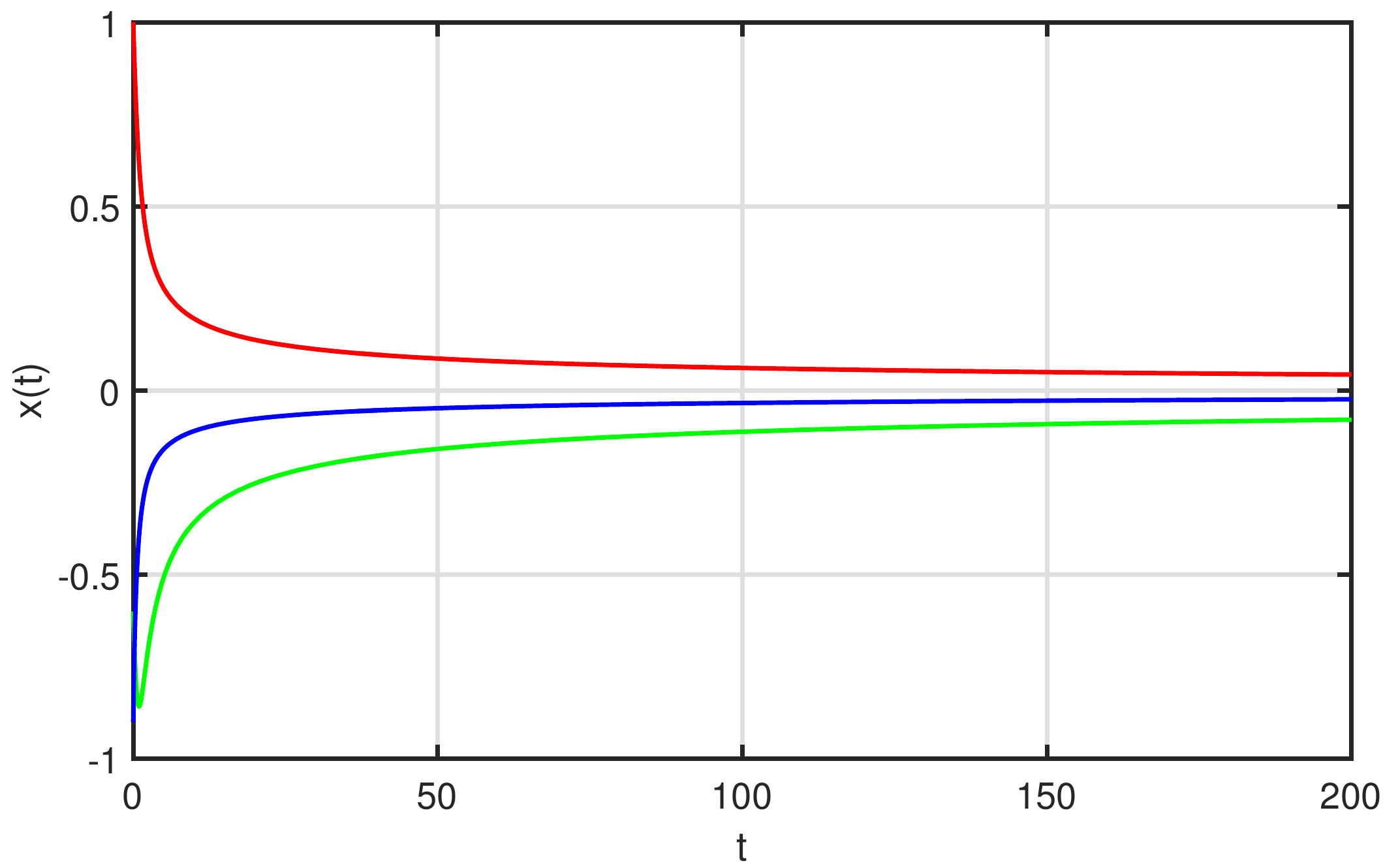}}
\caption{Stubborn neutrals scenario in Example~\ref{exm2}: (a) Evolution of each agent for all $x_i(0)>0$; (b) Evolution of each agent under the condition that $\exists i,j~\text{s.t.}~x_i(0)>0, x_j(0)<0$.}
\label{figure8}
\end{figure}

\begin{figure}[htbp]
\centering
\subfigure[]{\includegraphics[width=0.42\textwidth]{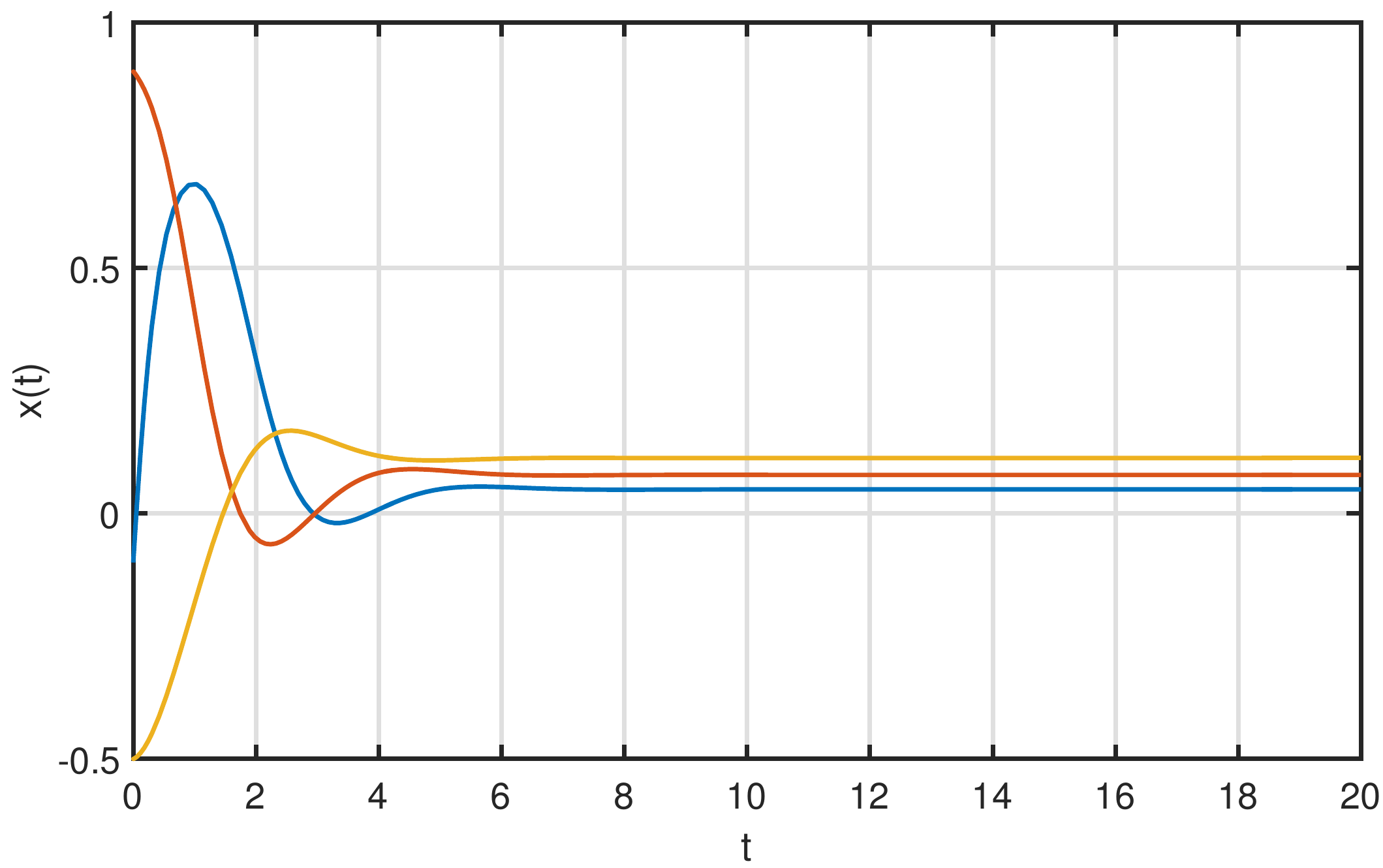}}
\subfigure[]{\includegraphics[width=0.42\textwidth]{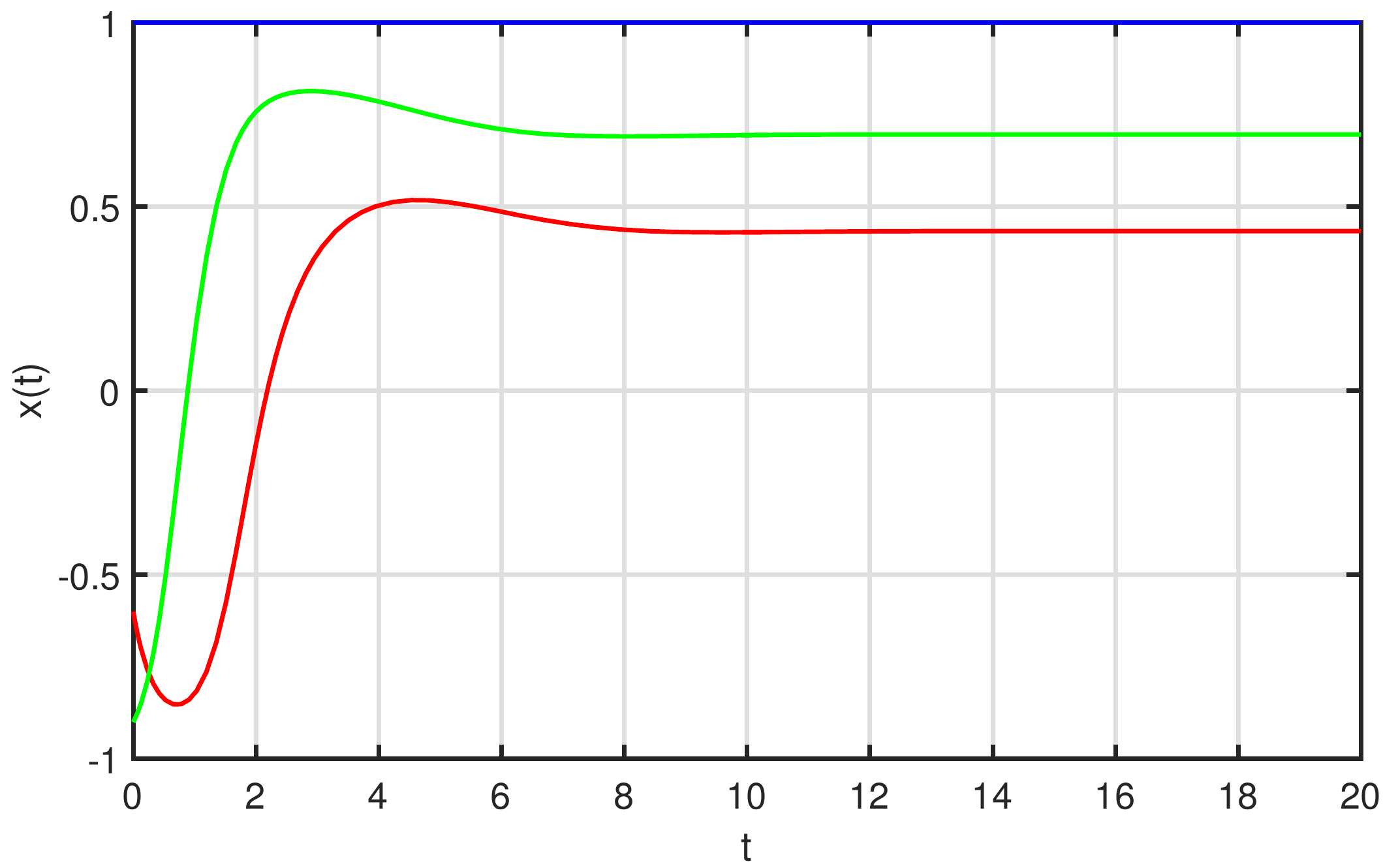}}
\subfigure[]{\includegraphics[width=0.42\textwidth]{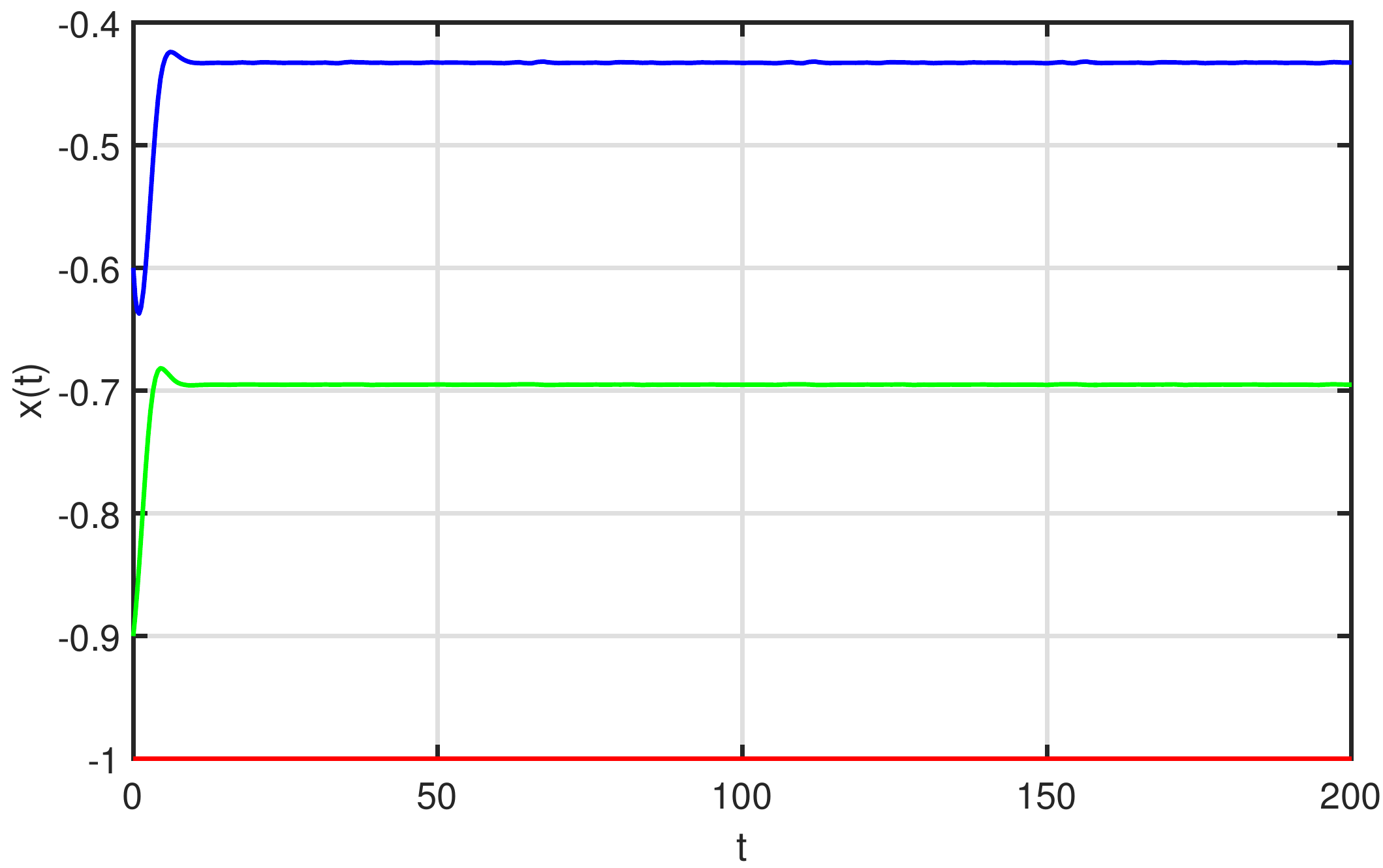}}
\subfigure[]{\includegraphics[width=0.42\textwidth]{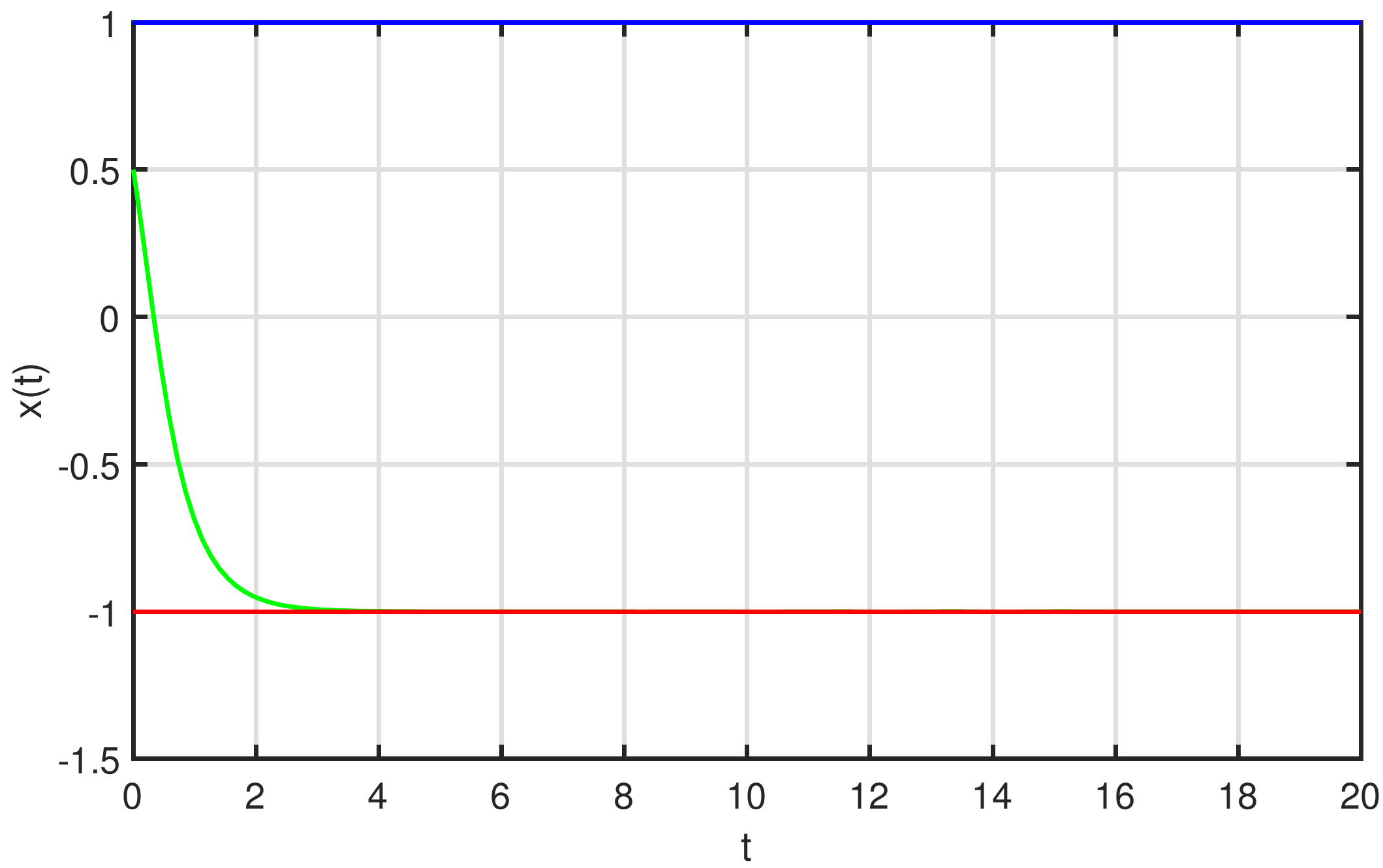}}
\caption{Stubborn extremists scenario in Example~\ref{exm2}: (a) Evolution of each agent for all $|x_i(0)|<1$; (b) Evolution of each agent under the condition that $x_3(0)=1, |x_i(0)|<1, \forall i\neq3$; (c) Evolution of each agent under the condition that $x_3(0)=-1, |x_i(0)|<1, \forall i\neq3$; (d) Evolution of each agent under the condition that $\exists i,j~\text{s.t.}~x_i(0)=1, x_j(0)=-1$.}
\label{figure9}
\end{figure}
\end{example}

\section{Conclusion}\label{section5}

A generalized nonlinear opinion dynamics with state-dependent susceptibility to persuasion and antagonistic interactions has been introduced in this paper. The introduced model has been comprehensively analyzed, particularly with a thorough investigation into three specializations of state-dependent susceptibility drawn from the existing literature. The main theoretical analysis tools adopted have been Perron-Frobenius property of eventually positive matrix and LaSalle invariance principle. It has been shown that all agents of the opinion network converge into the subspace spanned by the positive eigenvector of an eventually positive matrix. Generally speaking, all entries of the positive eigenvector are not the same unless the eventually positive matrix is a Laplacian matrix. Therefore, the obtained results can represent different levels of an opinion, which is a common phenomenon in social networks. Finally, in this paper, it has been assumed that all eigenvalues of the system matrix are non-positive. An interesting topic for future work will be to study the case that there does not exist any restriction on eigenvalues of a network matrix.

\appendices
\section{Signed Graphs}\label{appendix1}

Let $\mathcal{G}=(\mathcal{V},\mathcal{E},B)$ denote a weighted directed signed graph (signed digraph), where $\mathcal{V}=\{1,2,\cdots,N\}$ is the set of all nodes, $\mathcal{E}\in \mathcal{V}\times\mathcal{V}$ is the set of edges, and $B$ is an adjacency matrix which assigns real numbers to the edges. Let $\mathcal{G}(B)$ denote a graph with adjacency matrix $B$. An edge $(i,j)\in\mathcal{E}$ is directed from node $i$ to node $j$. The entry $b_{ij}>0\, (\mathrm{or }<0)$ is the weight corresponding to the edge $(i,j)$. Define $\mathcal{E}^{+}=\{(i,j)\mid b_{ij}>0\}$ and $\mathcal{E}^{-}=\{(i,j)\mid b_{ij}<0\}$. Then $\mathcal{E}=\mathcal{E}^{+}\cup\mathcal{E}^{-}$. For signed digraph $\mathcal{G}(B)$, let $L=C_r-B$ denote the Laplacian matrix of $\mathcal{G}(B)$, where $C_r=\diag\Big\{\sum\limits_{j=1}^N |b_{1j}|,\cdots,\sum\limits_{j=1}^N |b_{Nj}|\Big\}$. There exists an undirected graph $\mathcal{G}(B_u)$ with adjacency matrix $B_u=(B+B^{\rm T})/2$, and the corresponding Laplacian is $L_u=(L+L^{\rm T})/2=C_r-B_u$. Define $C_c=\diag\Big\{\sum\limits_{i=1}^N |b_{i1}|,\cdots,\sum\limits_{i=1}^N |b_{iN}|\Big\}$. It is obvious that $C_r=C_c$ when $B^{\rm T}=B$. A signed digraph $\mathcal{G}(B)$ is said to be weight balanced if $C_r=C_c$.

A directed path (length $l-1$) is a sequence of directed edges of the form $(i_1,i_2),(i_2,i_3),\cdots,(i_{l−1},i_l)$ with distinct nodes. A signed digraph has a spanning tree if there is a root node, which has directed paths to all other nodes. A signed digraph is said to be strongly connected if there is a directed path between any two distinct nodes. $\mathcal{G}(B)$ is strongly connected if and only if $B$ is irreducible.

\section{Eventually Positive Matrices and Perron–Frobenius Property}\label{appendix2}

\begin{definition}
A matrix $A\in\mathbb{R}^{n\times n}$ is said to possess Perron-Frobenius property if $\rho(A)$ is a positive eigenvalue of $A$ and the corresponding right eigenvector is nonnegative.
\end{definition}

\begin{definition}
A matrix $A\in\mathbb{R}^{n\times n}$ is said to possess the strong Perron-Frobenius property if $\rho(A)$ is a simple positive eigenvalue of $A$ and the corresponding right eigenvector is positive.
\end{definition}

\begin{definition}
A matrix $A\in\mathbb{R}^{n\times n}$ is said to be eventually positive (eventually nonnegative) if there exists a positive integer $k_0$ such that $A^k>0\, (A^k\geq 0), \forall k\geq k_0$.
\end{definition}

The following lemma presents some connection between the eventual positivity and strong Perron-Frobenius property.
\begin{lemma}[\cite{noutsos2006perron}]\label{lemma2}
For a matrix $A\in\mathbb{R}^{n\times n}$, the following statements are equivalent:
\begin{enumerate}[1)]
\item Both matrices $A$ and $A^{\rm T}$ possess the strong Perron-Frobenius property.
\item $A$ is an eventually positive matrix.
\item $A^{\rm T}$ is an eventually positive matrix.
\end{enumerate}
\end{lemma}

\begin{lemma}[\cite{altafini2015predictable}]\label{lemma3}
Consider an eventually positive matrix $A$, and denote a right eigenvector of $A$ by $v_r>0$. Then any eigenvector $v_1$ of $A$ with $v_1>0$ must be a multiple of $v_r$.
\end{lemma}

\section{Some Facts about Semidefinite Matrices}\label{appendix3}

\begin{lemma}[\cite{horn2012matrix}]\label{lemma4}
Let $A\in\mathbb{R}^{n\times n}$ be negative semidefinite and let $x\in\mathbb{R}^n$. Then $x^{\rm T}Ax=0$ if and only if $Ax=0$.
\end{lemma}

\begin{lemma}[\cite{horn2012matrix}]\label{lemma5}
Suppose that $A\in\mathbb{R}^{n\times n}$ and $H(A)=\frac{1}{2}(A+A^{\rm T})$ is negative semidefinite. Then
\begin{enumerate}[1)]
\item $\text{nullspace}A\subset \text{nullspace}H(A)$; $\text{nullspace}A^{\rm T}\subset \text{nullspace}H(A)$.

\item $\rank H(A)\leq \rank A$.

\item The following statements are equivalent:
\begin{enumerate}[(i)]
\item $A$ and $H(A)$ have the same null space.
\item $A^{\rm T}$ and $H(A)$ have the same null space.
\item $\rank H(A)=\rank A$.
\end{enumerate}
\end{enumerate}
\end{lemma}

\begin{lemma}\label{lemma6}
Suppose that $A\in\mathbb{R}^{n\times n}$ and $H(A)=\frac{1}{2}(A+A^{\rm T})$ is negative semidefinite. If $A$ is orthogonally diagonalizable, then $A$ and $H(A)$ have the same null space, and $A^{\rm T}$ and $H(A)$ have the same null space.
\end{lemma}

\begin{IEEEproof}\
If $A$ is orthogonally diagonalizable, then there exists a real orthogonal matrix $B$ such that $BAB^{\rm T}$ is a diagonal matrix. Meanwhile, $A^{\rm T}$ is orthogonally diagonalizable, and $BA^{\rm T}B^{\rm T}=[BAB^{\rm T}]^{\rm T}$. Hence, $\rank H(A)=\rank A$. By Lemma~\ref{lemma4}, $A$ and $H(A)$ have the same null space, and $A^{\rm T}$ and $H(A)$ have the same null space.
\end{IEEEproof}



\end{document}